\documentclass[12pt]{article}
\renewcommand{\baselinestretch}{1.0}
\usepackage{amsmath,amssymb,esint,amsthm,bbm}
\usepackage{verbatim}
\usepackage[OT2,T1]{fontenc}
\usepackage{mathrsfs}
\usepackage{graphicx}

\newcommand{\singlespace}{\renewcommand{\baselinestretch}{0.8}
\large\normalsize}

\def\down#1{_{{}_{\scriptstyle #1}}}
\def\vbar{\mathchoice{\vrule height6.3ptdepth-.5ptwidth.8pt\kern-.8pt}
   {\vrule height6.3ptdepth-.5ptwidth.8pt\kern-.8pt}
   {\vrule height4.1ptdepth-.35ptwidth.6pt\kern-.6pt}
   {\vrule height3.1ptdepth-.25ptwidth.5pt\kern-.5pt}}

\def\bbc#1#2{{\rm \mkern#2mu\vbar\mkern-#2mu#1}}
\def\bbb#1{{\rm I\mkern-3.5mu #1}}
\def\bba#1#2{{\rm #1\mkern-#2mu{\cal F}udge #1}}
\def\bb#1{{\count4=`#1 \advance\count4by-64 \ifcase\count4\or\bba A{11.5}\or
   \bbb B\or\bbc C{5}\or\bbb D\or\bbb E\or\bbb F \or\bbc G{5}\or\bbb H\or
   \bbb I\or\bbc J{3}\or\bbb K\or\bbb L \or\bbb M\or\bbb N\or\bbc O{5} \or
   \bbb P\or\bbc Q{5}\or\bbb R\or\bba S{8}\or\bba T{10.5}\or\bbc U{5}\or
   \bba V{12}
\or\bba W{16.5}\or\bba X{11}\or\bba Y{11.7}\or\bba Z{7.5}{\cal F}i}}

\newcounter{theorem}
\def \R {\mathbb{R}}

\def \k {{\kappa}}

\begin{document}

\def \Lip{{\rm Lip}}
\def \ds{\displaystyle}
 \newcommand{\bea}{\begin{eqnarray}}
 \newcommand{\vep}{{\varepsilon}}
    \newcommand{\X} {{\cal X}}
    \newcommand{\ep} {{\epsilon}}
    \newcommand{\g} {{\bf g}}
    \newcommand{\F} {{\cal F}}
    \newcommand{\J} {{\cal J}}
    \newcommand{\B} {{\cal B}}
    \newcommand{\CS}{{\cal{S}}}
    \newcommand{\CR}{{\cal{R}}}
    \newcommand{\Lq} {{\cal L}}
    \newcommand{\cend} {\end{center}}
    \newcommand{\elp}[2]{L^{#1}{(#2)}}
    \newcommand{\jac}{\mathcal{J}}
    \newcommand{\eea}{\end{eqnarray}}
    \newcommand{\rn}{\mathbb{R}^n}
    \newcommand{\fracc}{\displaystyle\frac}
    \newcommand{\intt}{\displaystyle{\int}}
\def\P{{\mathcal P}}
    \newcommand{\ddball}{\delta\mbox{-}d\mbox{-ball}}
    \newcommand{\ointt}{\displaystyle{\oint}}
    \newcommand{\disp}{\displaystyle}
    \newcommand{\limm}{\displaystyle\lim}
    \newcommand{\ra}{\rightarrow}
    \newcommand{\sumk} {\sum_{k=1}^{m_j}}
    \newcommand{\sumj} {\sum_{j=1}^\infty}
    \newcommand{\beginc} {\begin{center}}
    \newcommand{\ltwo} {L^2(\Omega)}
    \newcommand{\ltwoj} {{L^2 (\Omega_j)}}
    \newcommand{\cl}[2] {{\cal L}^{#1}_\mu(#2,{Q})}
    \newcommand{\wa}[2] {W^{1,{#1}}_{\nu,\mu}(#2,{Q})}
    \newcommand{\wb}[2] {W^{1,{#1}}_{\nu,\mu,0}(#2,{Q})}
    \newcommand{\wc}[2] {{\cal W}^{1,{#1}}_{\nu,\mu}(#2,{Q})}
    \newcommand{\we}[2] {{\cal W}^{1,{#1}}_{\nu,\mu,0}(#2,{Q})}
    \newcommand{\normal}[1]{\vec #1}
    \newcommand{\hf} {[f]_{\alpha \; ; \; x_0}}
    \newcommand{\di} {\partial}
    \newcommand{\scq} {\mathcal{Q}}
    \newcommand{\bg} {{\bf g}}
    \newcommand{\diam}{{\rm diam}}
    \newcommand{\ball}[3] {B_{#1}({#2}_{#3})}
    \newcommand{\sumjn} {\sum_{j=1}^N}
    \newcommand{\sumjm} {\sum_{j=1}^M}
    \newcommand{\scl} {\mathcal{L}}
    \newcommand{\rhu} {\rightharpoonup}
    \newcommand{\lhu} {\leftharpoonup}
    \newcommand{\limj} {\disp\lim_{j\rightarrow \infty}}
    \newcommand{\dist}{\textrm{dist}}

    \newtheorem{remark}[theorem]{{\bf Remark}}
    \newtheorem{defn}[theorem]{{\bf Definition}}
    \newtheorem{thm}[theorem]{\bf{Theorem}}
    \newtheorem{propn}[theorem]{\bf{Proposition}}
    \newtheorem{lemma}[theorem]{\bf{Lemma}}
    \newtheorem{fact}[theorem]{{\bf Fact}}
    \newtheorem{corollary}[theorem]{{\bf Corollary}}

\begin{center}{\bf\Large A Compact Embedding Theorem \\
for Generalized Sobolev Spaces}
\footnote{\noindent{\footnotesize 2010 Mathematics Subject
Classification: 46B50, 46E35, 35H20

{\noindent\footnotesize Key words and phrases: compact
embedding, Sobolev spaces, degenerate quadratic forms}}} \end{center}
\begin{center}
    by Seng-Kee Chua, Scott Rodney and Richard L. Wheeden
\end{center}
\vspace{0.3cm}
{\footnotesize{{\bf Abstract:}  We give an elementary proof
    of a compact embedding theorem in abstract Sobolev spaces. The
    result is first presented in a general context and later
    specialized to the case of degenerate Sobolev spaces defined
    with respect to nonnegative quadratic forms on
    $\mathbb{R}^n$. Although our primary interest concerns degenerate
    quadratic forms, our result also applies to nondegenerate cases,
    and we consider several such applications, including the classical
    Rellich-Kondrachov compact embedding theorem and results for the
    class of $s$-John domains in $\mathbb{R}^n$, the latter for
    weights equal to powers of the distance to the boundary. We also
    derive a compactness result for Lebesgue spaces on quasimetric
    spaces unrelated to $\mathbb{R}^n$ and possibly without any notion
    of gradient.}}

\section{The General Theorem}

The main goal of this paper is to generalize the classical
Rellich-Kondrachov theorem concerning compact embedding of Sobolev spaces
into Lebesgue spaces. Our principal result applies not only to the
classical Sobolev spaces on open sets $\Omega\subset \rn$ but
also allows us to treat the degenerate Sobolev spaces defined in
\cite{SW2}, and to obtain compact embedding of them into various
$L^q(\Omega)$ spaces. These degenerate Sobolev spaces are associated
with quadratic forms $Q(x,\xi) = \xi' Q(x)\xi$, $x\in \Omega, \xi \in
\rn$, which are nonnegative but may vanish identically in $\xi$ for
some values of $x$. Such quadratic forms and Sobolev spaces arise
naturally in the study of existence and regularity of weak solutions
of some second order subelliptic linear/quasilinear partial
differential equations; see, e.g., \cite[2]{SW1}, \cite{R1},
\cite{MRW}, \cite{RSW}.

The Rellich-Kondrachov theorem is frequently used to study the
existence of solutions to elliptic equations, a famous example being
subcritical and critical Yamabe equations, resulting in the solution
of Yamabe's problem; see \cite{Y}, \cite{T}, \cite{A},
\cite{S}. Further applications lie in proving the existence of weak
solutions to Dirichlet problems for elliptic equations with rough
boundary data and coefficients; see \cite{GT}. In a sequel to this
paper, we will apply our compact embedding results to study the
existence of solutions for some classes of degenerate equations.

In this section, we will state and prove our most general compact
embedding results. In Sections 2 and 3, we study some applications to
classical and degenerate Sobolev spaces, respectively. In Section 4,
more general results in quasimetric spaces are studied.

We begin by listing some useful notation.  Let $w$ be
a measure on a $\sigma$-algebra $\Sigma$ of subsets of a set
$\Omega$, with $\Omega \in \Sigma$. For $0< p\le \infty$, let
$L^p_w(\Omega)$ denote the class of real-valued measurable functions
$f$ satisfying $||f||_{L^p_w(\Omega)} <\infty$, where
$||f||_{L^p_w(\Omega)} =  \Big(\int_\Omega |f|^pdw
\Big)^{1/p}$ if $p <\infty$ and $||f||_{L^\infty_w(\Omega)} =
\text{ess sup}_\Omega\, |f|$, the essential supremum being taken with
respect to $w$-measure.  When dealing with generic functions in
$L^p_w(\Omega)$, we will not distinguish between functions which are
equal a.e.-$w$. For $E\in \Sigma$, $w(E)$ denotes the
$w$-measure of $E$, and if $0<w(E)<\infty$ then $f_{E,w}$ denotes the
$w$-average of $f$ over $E$: $f_{E,w} = \frac{1}{w(E)}\int_E fdw$.
Throughout the paper, positive constants will be denoted by $C$ or $c$
and their dependence on important parameters will be indicated.

For $k\in\mathbb{N}$, let $\mathscr{X}(\Omega)$ be a normed linear
space of measurable $\mathbb{R}^k$-valued functions $\g$ defined on
$\Omega$ with norm $||\g||_{\mathscr{X}(\Omega)}$.  We assume that
there is a subset $\Sigma_0\subset \Sigma$ so that
$(\mathscr{X}(\Omega), \Sigma_0)$ satisfies the following
properties:\\

\noindent{(A)} For any $\g\in\mathscr{X}(\Omega)$ and $F\in
\Sigma_0$, the function $\g\chi_F\in \mathscr{X}(\Omega)$, where
$\chi_F$ denotes the characteristic function of $F$.\\

\noindent{($B_p$)} There are constants $C_1, C_2, p$ satisfying $1
\le C_1, C_2, p < \infty$ so that if $\{F_\ell\}$ is
a finite collection of sets in $\Sigma_0$ with $\disp\sum_{\ell}
\chi_{F_\ell}(x) \leq C_1$ for all $x\in\Omega$, then
\bea
\disp\sum_{\ell} ||\g\chi_{F_\ell}||^{p}_{\mathscr{X}(\Omega)} \leq
C_2  ||\g||^{p}_{\mathscr{X}(\Omega)}\quad \text{for all $\g\in
\mathscr{X}(\Omega)$.} \nonumber
\eea

For $1\le N \le \infty$, we will often consider the product space
$L^N_w(\Omega)\times \mathscr{X}(\Omega)$.  This is a normed linear
space with norm
\bea\label{bspace} ||(f,\g)||_{L^N_w(\Omega)\times \mathscr{X}(\Omega)} =
||f||_{L^N_w(\Omega)} + ||\g||_{\mathscr{X}(\Omega)}.
\eea
A set ${\cal S}\subset L^N_w(\Omega)\times \mathscr{X}(\Omega)$ will
be called a {\it bounded set in} $L^N_w(\Omega)\times
\mathscr{X}(\Omega)$ if
$$
\sup_{(f,\g)\in {\cal S}} ||(f,\g)||_{L^N_w(\Omega)\times
\mathscr{X}(\Omega)} < \infty.
$$
Projection maps such as the one defined by
\bea\label{projection}
\pi:(f,\g) \ra f, \quad (f,\g) \in L^N_w(\Omega)\times
\mathscr{X}(\Omega),
\eea
will play a role in our results. If $w(\Omega) < \infty$, then
$\pi(L^N_w(\Omega)\times \mathscr{X}(\Omega))\subset  L^q_w(\Omega)$
if $1\le q \le N.$

\begin{thm}\label{general} Let $w$ be a finite measure on a
$\sigma$-algebra $\Sigma$ of subsets of a set $\Omega$, with $\Omega
\in \Sigma$. Let $1 \le p < \infty$, $1<N \le \infty$,
$\mathscr{X}(\Omega)$ be a normed linear space satisfying properties
(A) and ($B_p$) relative to a collection $\Sigma_0 \subset \Sigma$,
and let ${\cal S}$ be a bounded set in $L^N_w(\Omega)\times
\mathscr{X}(\Omega)$. 

Suppose that ${\cal S}$ satisfies the following: given $\epsilon>0$,
there are a finite number of pairs $\{E_\ell,F_\ell\}_{\ell=1}^J$ with
$E_\ell \in \Sigma$ and $F_\ell\in \Sigma_0$ (the pairs and $J$ may
depend on $\epsilon$) such that

\noindent (i) $w\big(\Omega\setminus\cup_{\ell} E_\ell\big) <\epsilon$
and $w(E_\ell) >0$;

\noindent (ii) $\{F_\ell\}$ has bounded overlaps independent of
$\epsilon$ with the same overlap constant as in ($B_p$), i.e.,
\bea\label{bddov} \disp\sum_{\ell=1}^J \chi_{F_\ell}(x) \leq C_1, \quad
x\in \Omega,
\eea
for $C_1$ as in ($B_p$);

\noindent (iii) for every $(f,\g)\in {\cal S}$, the local
Poincar\'e-type inequality
\bea\label{pc} ||f-f_{E_\ell,w}||_{L^{p}_w(E_\ell)} \leq \epsilon\,
||\g\chi_{F_\ell}||_{\mathscr{X}(\Omega)}
\eea
holds for each $(E_\ell, F_\ell)$.

Let $\hat{\cal S}$ be the set defined by
\bea\label{hatS}
\hat{\cal S}= \left\{f\in L^N_w(\Omega):\, \text{there exists
$\{(f^j,\g^j)\}_{j=1}^\infty \subset {\cal S}$ with $f^j \to f$
a.e.-$w$}\right\}.
\eea
Then $\hat{\cal S}$ is compactly embedded in $L^q_w(\Omega)$ if $1\le q
<N$ in the sense that for every sequence $\{f_k\} \subset \hat{\cal
S}$, there is a single subsequence $\{f_{k_i}\}$ and a function $f\in
L^N_w(\Omega)$ such that $f_{k_i}\to f$ pointwise a.e.-$w$ in $\Omega$
and in $L^q_w(\Omega)$ norm for $1\leq q<N$.
\end{thm}

Before proceeding with the proof of Theorem \ref{general}, we make
several simple observations. First, in the definition of $\hat{\cal
S}$, the property that $f\in L^N_w(\Omega)$ follows by Fatou's lemma
since the associated functions $f^j$ are bounded in $L^N_w(\Omega)$,
as ${\cal S}$ is bounded in $L^N_w(\Omega) \times\mathscr{X}(\Omega)$
by hypothesis. Fatou's lemma also shows that $\hat{\cal S}$ is a
bounded set in $L^N_w(\Omega)$. Moreover, since $N>1$, if $\{f^j\}$ is
bounded in $L^N_w(\Omega)$ and $f^j\to f$ a.e.-$w$, then $(f^j)_{E,w}
\to f_{E,w}$ for all $E\in \Sigma$; in fact, in this situation, by
using Egorov's theorem, we have $\int_\Omega f^j \varphi dw \to \int_\Omega
f\varphi dw$ for all $\varphi \in L^{N'}_w(\Omega), 1/N +1/N'=1$.

Next, while the hypothesis $w(E_\ell)>0$ in assumption (i)
ensures that the averages $f_{E_\ell,w}$ in (\ref{pc}) are
well-defined, it is not needed since we can discard any pair
$E_\ell, F_\ell$ with $w(E_\ell) =0$ without affecting the inequality
$w(\Omega\setminus \cup E_\ell)<\epsilon$ or (\ref{bddov}) and
(\ref{pc}).

Finally, since $\hat{\cal S}$ contains the first component
$f$ of any pair $(f,\g) \in {\cal S}$, a simple corollary of Theorem
\ref{general} is that the projection $\pi$ defined in
(\ref{projection}) is a compact mapping of ${\cal S}$ into
$L^q_w(\Omega)$, $1\le q <N$, in the sense that for every
sequence $\{(f_k,\g_k)\}\subset {\cal S}$, there is a subsequence
$\{f_{k_i}\}$ and a function $f\in L^N_w(\Omega)$ such that
$f_{k_i}\to f$ pointwise a.e.-$w$ in $\Omega$ and in $L^q_w(\Omega)$
norm for $1\leq q<N$.  \\

\noindent {\bf Proof:}  Let ${\cal S}$ satisfy the hypotheses and
suppose $\{f_k\}_{k\in\mathbb{N}} \subset \hat{\cal S}$. For each
$f_k$, use the definition of $\hat{\cal S}$ to choose a sequence
$\{(f_k^j, \g_k^j)\}_j \subset {\cal S}$ with $f_k^j \to f_k$
a.e.-$w$ as $j\to \infty$. Since ${\cal S}$ is bounded in
$L^N_w(\Omega)\times \mathscr{X}(\Omega)$, there is $M\in (0,\infty)$
so that $||(f_k^j,\g_k^j)||_{L^N_w(\Omega)\times \mathscr{X}(\Omega)}
\leq M$ for all $k$ and $j$. Also, as noted above, $\{f_k\}$ is
bounded in $L^N_w(\Omega)$ norm; in fact $||f_k||_{L^N_w(\Omega)} \le
M$ for the same constant $M$ and all $k$.

Since $\{f_k\}$ is bounded in $L^N_w(\Omega)$, then if
$1<N<\infty$, it has a weakly convergent subsequence, while if
$N=\infty$, it has a subsequence which converges in the weak-star
topology. In either case, we relabel the subsequence
as $\{f_k\}$ to preserve the index. Fix $\epsilon >0$ and let
$\{E_\ell, F_\ell\}_{\ell =1}^J$ satisfy the hypotheses of the theorem
relative to $\epsilon$. Setting $\Omega^\epsilon = \cup E_\ell$, we
have by assumption (i) that
\bea \label{Oj}
w(\Omega\setminus\Omega^\epsilon) < \epsilon.
\eea

Let us show that there is a positive constant $C$ independent of
$\epsilon$ so that
\bea\label{diff}
\disp\sum_\ell ||f_k - (f_k)_{E_\ell,w}||_{L^p_w(E_\ell)}^p \le
C\epsilon^p \quad\mbox{for all $k$.}
\eea
Fix $k$ and let $\Delta$ denote the expression on the left side of
(\ref{diff}). Since $f_k^j - (f_k^j)_{E_\ell,w} \to f_k -
(f_k)_{E_\ell,w}$ a.e.-$w$ as $j \to \infty$, Fatou's lemma gives
\bea
\Delta \le \disp\sum_\ell \liminf_{j\to\infty} ||f_k^j -
(f_k^j)_{E_\ell,w}||_{L^p_w(E_\ell)}^p.\nonumber
\eea
Consequently, by using the Poincar\'e inequality (\ref{pc}) for ${\cal
S}$ and superadditivity of $\liminf$, we obtain
\bea
\Delta \le  \liminf_{j\to \infty}  \disp\sum_\ell \epsilon^p ||\g_k^j
\chi_{F_\ell} ||^p_{\mathscr{X}(\Omega)}.\nonumber
\eea
By (\ref{bddov}), the sets $F_\ell$ have finite overlaps uniformly in
$\epsilon$, with the same overlap constant $C_1$ as in property ($B_p$) of
$\mathscr{X}(\Omega)$. Hence, by property ($B_p$) applied to the last
expression together with boundedness of ${\cal S}$,
\bea
\Delta \le C_2 \epsilon^p \liminf_{j\to \infty}
||\g_k^j||_{\mathscr{X}(\Omega)}^p  \le C_2M^p\epsilon^p.\nonumber
\eea
This proves (\ref{diff}) with $C= C_2M^p$.

Next note that
\[ \int_{\Omega^\epsilon}|f_m-f_k|^{p}dw \le\,
\sum_{\ell}\int_{E_\ell}|f_m-f_k|^{p}dw
\]
\[
\le 2^{p-1}\Big(\sum_{\ell} \int_{E_\ell}|f_m-f_k
-(f_m-f_k)_{E_\ell,w}|^{p}dw +
\disp\sum_{\ell}|(f_m-f_k)_{E_\ell,w}|^{p}w(E_\ell)\Big)
\]
\bea \label{IandII}
= 2^{p-1}(I + II).
\eea
We will estimate I and II separately. We have
\bea
I \le  2^{p-1}\left(\disp\sum_{\ell} ||f_m -
(f_m)_{E_\ell,w}||_{L^p_w(E_\ell)}^p + \disp\sum_\ell
||f_k - (f_k)_{E_\ell,w}||_{L^p_w(E_\ell)}^p\right)\nonumber
\eea
\bea \label{I}
\leq 2^{p-1}\left(C\epsilon^{p} + C\epsilon^p \right) = 2^{p}C
\epsilon^{p}
\eea
by (\ref{diff}). To estimate $II$, first note that
\[
II = \disp\sum_{\ell=1}^J |(f_m-f_k)_{E_\ell,w}|^{p}w(E_\ell) =
\disp\sum_{\ell=1}^J \frac{1}{w(E_l)^{p-1}}\Big|\int_\Omega
(f_m-f_k)\chi_{E_\ell}dw\Big|^{p}.
\]
Since $w(\Omega)<\infty$, each characteristic function
$\chi_{E_\ell}\in L^{N'}_w(\Omega)$, $1/N + 1/N' =1$ (with $N'=1$ if
$N=\infty$).  As $\{f_k\}$ converges weakly in $L^N_w(\Omega)$ when
$1<N<\infty$, or converges in the weak-star sense when $N=\infty$,
then for $m,k$ sufficiently large depending on $\epsilon$, and for all
$1\leq \ell \leq J$,
\[
\frac{1}{w(E_l)^{p-1}}\Big| \int_{\Omega} (f_m-f_k)
\chi_{E_\ell}dw\Big|^{p}  \leq  \fracc{\epsilon^{p}}{J}.
\]
Thus $II\leq \epsilon^{p}$ for $m,k$ sufficiently large
depending on $\epsilon$.  Combining this estimate with (\ref{IandII})
and (\ref{I}) shows that
\bea \label{COj} ||f_m-f_k||_{L^{p}_w(\Omega^\epsilon)}< C\epsilon
\eea
for $m,k$ sufficiently large and $C =C(M, C_2)$.

Let us now show that $\{f_k\}$ is a Cauchy sequence in
$L^{1}_w(\Omega)$. For $m, k$ as in (\ref{COj}), H\"older's inequality
and the fact that $||f_k||_{L^N_w(\Omega)} \leq M$ for all $k$ yield
\bea ||f_m-f_k||_{L^{1}_w(\Omega)} &\le & ||f_m-f_k
||_{L^{1}_w(\Omega^\epsilon)} + ||f_m-f_k||_{L^{1}_w(\Omega
\setminus \Omega^\epsilon)}\nonumber\\ &\leq& ||f_m-
f_k||_{L^{p}_w(\Omega^\epsilon)} w(\Omega^\epsilon)^{\frac{1}{p'}}  +
||f_m-f_k||_{L^{N}_w(\Omega \setminus\Omega^\epsilon)}
w(\Omega\setminus \Omega^\epsilon)^{\frac{1}{N'}}\nonumber\\  &<&
C\epsilon w(\Omega^\epsilon)^{\frac{1}{p'}} + 2Mw(\Omega\setminus
\Omega^\epsilon)^{\frac{1}{N'}}\nonumber\\ &<& C\epsilon
w(\Omega)^\frac{1}{p'}  +2M\epsilon^{\frac{1}{N'}}  \quad \mbox{by
(\ref{Oj})}.\nonumber
\eea
Since $N'<\infty$, it follows that $\{f_k\}$ is Cauchy in
$L^{1}_w(\Omega)$. Hence it has a subsequence (again denoted by
$\{f_k\}$) that converges in $L^1_w(\Omega)$ and pointwise a.e.-$w$ in
$\Omega$ to a function $f\in L^1_w(\Omega)$. If $N=\infty$,
$\{f_k\}$ is bounded in $L^\infty_w(\Omega)$ by hypothesis, so its
pointwise limit $f \in L^\infty_w(\Omega)$. If $N<\infty$, since
$\{f_k\}$ is bounded in $L^N_w(\Omega)$, Fatou's Lemma implies that
$f\in L^N_w(\Omega)$. This completes the proof in case $q=1$.

For general $q$, we will use the same subsequence $\{f_k\}$ as
above. Thus we only need to show that $\{f_k\}$ converges in
$L^q_w(\Omega)$ for $1<q<N$. We will use H\"older's inequality.
Given $q\in (1,N)$, choose $\lambda\in(0,1)$, namely $\lambda =
\big(\frac{1}{q}-\frac{1}{N}\big)/\big(1 -\frac{1}{N}\big)$, hence
$\lambda = 1/q$ if $N=\infty$, so that
\bea \label{interp} ||f_m-f_k||_{L^q_w(\Omega)} \leq
||f_m-f_k||_{L^{1}_w(\Omega)}^\lambda
||f_m-f_k||_{L^N_w(\Omega)}^{1-\lambda}.
\eea
As before, $||f_k||_{L^N_w(\Omega)} \le M$, and therefore
$||f_m-f_k||_{L^N_w(\Omega)}^{1-\lambda}\leq (2M)^{1-\lambda}$, giving
by (\ref{interp}) that $\{f_k\}$ is Cauchy in $L^q_w(\Omega)$ as it is
Cauchy in $L^{1}_w(\Omega)$.  This completes the proof of Theorem
\ref{general}. $\Box$

A compact embedding result is also proved in \cite[Theorem
3.4]{FSSC} by using Poincar\'e type estimates. However, Theorem
\ref{general} applies to situations not considered in \cite{FSSC}
since it is not restricted to the context of Lipschitz vector fields
in $\mathbb{R}^n$. Other abstract compact embedding results can be
found in \cite[Theorem 4]{HK1} and \cite[Theorem 8.1]{HK2}, including
a version (see \cite[Theorem 5]{HK1}) for weighted Sobolev spaces with
nonzero continuous weights, and a version in \cite{HK2} for metric
spaces with a single doubling measure. The proof in \cite{HK1} assumes
prior knowledge of the classical Rellich-Kondrachov compactness
theorem (see e.g. \cite[Theorem 7.22(i)]{GT} and below).\\

By making minor changes in the proof of Theorem \ref{general}, we can
obtain a sufficient condition for a bounded set in $L^N_w(\Omega)$ to
be precompact in $L^q_w(\Omega)$, $1\le q< N$, without mentioning
the sets $\{F_\ell\}$, the space $\mathscr{X}(\Omega)$, properties (A)
and ($B_p$), or conditions (\ref{bddov}) and (\ref{pc}). We state this
result in the next theorem. An application is given in \S 4.

\begin{thm}\label{absgeneral} Let $w$ be a finite measure on a
$\sigma$-algebra $\Sigma$ of subsets of a set $\Omega$, with $\Omega
\in \Sigma$. Let $1 \le p <\infty$, $1< N \le \infty$ and $\P$ be a
bounded subset of $L^{N}_w(\Omega)$. Suppose there is a positive
constant $C$ so that for every $\epsilon>0$, there are a finite number
of sets $E_\ell \in \Sigma$  with

\noindent (i) $w\big(\Omega\setminus\cup_{\ell} E_\ell\big) <\epsilon$
and $w(E_\ell) >0$;

\noindent (ii)  for every $f \in \P$,
\begin{equation}\label{abspoincare}
\sum_{\ell} ||f-f_{E_\ell,w}||^p_{L^{p}_w(E_\ell)}\le C\epsilon^p.
\end{equation}

\noindent Let
$$\hat{\P}=\{ f\in L^N_w(\Omega): \ \mbox{there exists }
\{f^j\}\subset \P \ \mbox{with $ f^j \to f \ a.e.$-$w$ } \}.
$$
Then for every sequence $\{f_k\}\subset \hat{\P}$, there is a single
subsequence $\{f_{k_i}\}$ and a function $f\in L^N_w(\Omega)$ such
that $f_{k_i}\to f$ pointwise a.e.-$w$ in $\Omega$ and in
$L^q_w(\Omega)$ norm for $1\leq q<N$.
\end{thm}

\begin{remark}\label{closureS}
\begin{enumerate}

\item Given $\epsilon >0,$ let $\{E_\ell\}$ satisfy hypothesis (i)
of Theorem \ref{absgeneral}. Hypothesis (ii) of Theorem
\ref{absgeneral} is clearly true for $\{E_\ell\}$ if for every $f\in
\P$, there are nonnegative constants $\{a_\ell\}$ such that
\bea\label{pc2} ||f-f_{E_\ell,w}||_{L^{p}_w(E_\ell)} \leq \epsilon\
a_\ell
\eea
and
\begin{equation}\label{sumcond}
\sum a_\ell^p \le C
\end{equation}
with $C$ independent of $f, \epsilon$. The constants $\{a_\ell\}$ may
vary with $f$ and $\epsilon$.

\item Theorem \ref{general} is a corollary of Theorem
\ref{absgeneral}. To see why, suppose that the hypothesis of Theorem
\ref{general} holds. Define $\P$ by $\P = \pi({\cal S}) = \{f:
(f,\g)\in {\cal S}\}$. Let $\epsilon>0$ and choose $\{(E_\ell,
F_\ell)\}$ as in Theorem \ref{general}. Given $f\in \P$, choose any
$\g$ such that  $(f,\g) \in {\cal S}$ and set $a_\ell =
||\g\chi_{F_\ell}||_{\mathscr{X}(\Omega)}$ for all $\ell$.  Then
(\ref{pc}), (\ref{bddov}) and property ($B_p$) of $\mathscr{X}(\Omega)$
imply (\ref{pc2}) and (\ref{sumcond}). The preceding remark shows
that the hypothesis of Theorem \ref{absgeneral} holds. The conclusion
of Theorem \ref{general} now follows from Theorem \ref{absgeneral}.
\end{enumerate}
\end{remark}

{\bf Proof of Theorem \ref{absgeneral}: } Theorem \ref{absgeneral} can
be proved by checking through the proof of Theorem \ref{general}. In
fact, the nature of hypothesis (\ref{abspoincare}) allows
simplification of the proof. First recall that if $f^j\to f$ a.e.-$w$
and $\{f^j\}$ is bounded in $L^N_w(\Omega)$, then $(f^j)_{E,w}\to
f_{E, w}$ for every $E\in \Sigma$. Therefore, by the definition of
$\hat{\P}$ and Fatou's lemma, the truth of (\ref{abspoincare}) for
all $f\in \P$ implies its truth for all $f \in \hat{\P}$.  Given a
sequence $\{f_k\}$ in $\hat{\P}$, we follow the proof of Theorem
\ref{general} but no longer need to introduce the $\{f_k^j\}$ or
prove (\ref{diff}) since (\ref{diff}) now follows from the fact
that (\ref{abspoincare}) holds for $\hat{\P}$. Further details are
left to the reader. \qed

We close this section by listing an alternate version of Theorem
\ref{general} that we will use in \S 3.4 when we consider local
results.

\begin{thm}\label{generallocal} Let $w$ be a measure (not necessarily
finite) on a $\sigma$-algebra $\Sigma$ of subsets of a set $\Omega$,
with $\Omega \in \Sigma$. Let $1 \le p < \infty$, $1<N \le \infty$,
$\mathscr{X}(\Omega)$ be a normed linear space satisfying properties
(A) and ($B_p$) relative to a set $\Sigma_0 \subset \Sigma$, and let
${\cal S}$ be a collection of pairs $(f, \bf{g})$ such that $f$ is
$\Sigma$-measurable and ${\bf g} \in \mathscr{X}(\Omega)$.

Suppose that ${\cal S}$ satisfies the following conditions relative to
a fixed set $\Omega'\in \Sigma$ (in particular $\Omega'\subset
\Omega$): for each $\epsilon =\epsilon_j= 1/j$ with $j\in \mathbb{N}$,
there are a finite number of pairs $\{E_\ell^\epsilon,
F_\ell^\epsilon\}_{\ell}$ with $E_\ell^\epsilon \in \Sigma$ and
$F_\ell^\epsilon \in \Sigma_0$ such that

\noindent (i) $w(\Omega' \setminus \cup_\ell E_\ell^\epsilon)=0$
and $0<w(E_\ell^\epsilon) <\infty$;

\noindent (ii) $\{F_\ell^\epsilon\}_\ell$ has bounded overlaps
independent of $\epsilon$ with the same overlap constant as in
($B_p$), i.e.,
\bea\nonumber \disp\sum_\ell \chi_{F_\ell^\epsilon}(x) \leq C_1, \quad
x\in \Omega,
\eea
for $C_1$ as in ($B_p$);

\noindent (iii) for every $(f,\g)\in {\cal S}$, the local
Poincar\'e-type inequality
\bea\nonumber ||f-f_{E_\ell^\epsilon,w}||_{L^{p}_w(E_\ell^\epsilon)}
\leq \epsilon\, ||\g\chi_{F_\ell^\epsilon}||_{\mathscr{X}(\Omega)}
\eea
holds for each $(E_\ell^\epsilon, F_\ell^\epsilon)$.

Then for every sequence $\{(f_k, \g_k)\}$ in ${\cal S}$ with
\bea\label{bddlocal}
\sup_k\left[||f_k||_{L^N_w(\cup_{\ell, j}E_\ell^{1/j})} +
||\g_k||_{\mathscr{X}(\Omega)} \right] < \infty,
\eea
there is a subsequence $\{f_{k_i}\}$ of $\{f_k\}$ and a function $f\in
L^N_w(\Omega')$ such that $f_{k_i}\to f$ pointwise a.e.-$w$ in $\Omega'$
and in $L^q_w(\Omega')$ norm for $1\leq q\le p$.  If $p<N$, then also
$f_{k_i}\to f$ in $L^q_w(\Omega')$ norm for $1\leq q <N$.
\end{thm}

The principal difference between the assumptions in Theorems
\ref{general} and \ref{generallocal} occurs in hypothesis (i).  When
we apply Theorem \ref{generallocal} in \S 3.4, the sets
$\{E_\ell^\epsilon\}$ will satisfy $\Omega' \subset
\cup_\ell E_\ell^\epsilon$ for each $\epsilon$, and consequently the
condition in hypothesis (i) that $w(\Omega'\setminus \cup_\ell
E^\epsilon_\ell)=0$ for each $\epsilon$ will be automatically true.
Unlike Theorem \ref{general}, the value of $q$ in Theorem
\ref{generallocal} is always allowed to equal $p$.  Although
$w(\Omega)$ is not assumed to be finite in Theorem \ref{generallocal},
$w(\Omega')<\infty$ is true due to hypothesis (i) and the fact that
the number of $E_\ell^\epsilon$ is finite for each $\epsilon$. As in
Theorem \ref{general}, the hypothesis $w(E_\ell^\epsilon)>0$ is
dispensible.

{\bf Proof of Theorem \ref{generallocal}: }
The proof is like that of Theorem \ref{general}, with minor changes
and some simplifications. We work directly with the pairs $(f_k,\g_k)$
without considering approximations $(f_k^j, \g_k^j)$. Due to the form
of assumption (i) in Theorem \ref{generallocal}, neither the set
$\Omega^\epsilon$ nor estimate (\ref{Oj}) is now needed. Since
$w(\Omega' \setminus \cup_\ell E_\ell^\epsilon)=0$ for each
$\epsilon =1/j$, we can replace $\Omega^\epsilon$ by $\Omega'$ in the
proof, obtaining the estimate
\bea\label{newCOj}
||f_m-f_k||_{L^{p}_w(\Omega')}< C\epsilon
\eea
as an analogue of (\ref{COj}). In deriving (\ref{newCOj}), the
weak and weak-star arguments are guaranteed since by (\ref{bddlocal}),
\[
\sup_k ||f_k||_{L^N_w(\cup_{\ell, j}E_\ell^{1/j})} <\infty.
\]
The main change in the proof comes by observing that the entire
argument formerly used to show that $\{f_k\}$ is Cauchy in
$L^1_w(\Omega)$ is no longer needed. In fact, (\ref{newCOj}) proves
that $\{f_k\}$ is Cauchy in $L^p_w(\Omega')$, and therefore it is also
Cauchy in $L^q_w(\Omega')$ if $1\le q \le p$ since $w(\Omega')
<\infty$. The first conclusion in Theorem \ref{generallocal} then
follows. To prove the second one, assuming that $p,q <N$, we use an
analogue of (\ref{interp}) with $\Omega'$ in place of $\Omega$ and the
same choice of $\lambda$, namely,
\[
||f_m-f_k||_{L^q_w(\Omega')} \leq
||f_m-f_k||_{L^1_w(\Omega')}^\lambda
||f_m-f_k||_{L^N_w(\Omega')}^{1-\lambda}.
\]
The desired conclusion then follows as before since we have already
shown that the first factor on the right side tends to $0$.

\section{Applications in the Nondegenerate Case}
\setcounter{equation}{0}

\setcounter{theorem}{0}

Roughly speaking, a consequence of Theorem \ref{general} is that
a set of functions which is bounded in $L^N_w(\Omega)$ is precompact
in $L^q_w(\Omega)$ for $1\le q<N$ if the gradients of the functions
are bounded in an appropriate norm, and a {\it local} Poincar\'e
inequality holds for them. The requirement of boundedness in
$L^N_w(\Omega)$ will be fulfilled if, for example, the functions
satisfy a {\it global} Poincar\'e or Sobolev estimate with exponent
$N$ on the left-hand side. In order to illustrate this principle more
precisely, we first consider the classical gradient operator and
functions on $\mathbb{R}^n$ with the standard Euclidean metric. We
include a simple way to see that the Rellich-Kondrachov compactness
theorem follows from our results. Our derivation of this fact is
different from those in \cite{AF} and \cite{GT}; in particular, it
avoids using the Arzel\'a-Ascoli theorem and regularization of
functions by convolution. We also list compactness results for the
special class of $s$-John domains in $\mathbb{R}^n$. In \cite{HK1},
the authors mention that such results follow from their development
without giving specific statements. See also \cite[Theorem
8.1]{HK2}. We list results for degenerate quadratic forms and vector
fields in Section 3.

We begin by proving a compact embedding result for some Sobolev spaces
involving two measures. Let $w$ be a measure on the Borel subsets of a
fixed open set $\Omega\subset \mathbb{R}^n$, and let $\mu$ be a
measure on the $\sigma$-algebra of Lebesgue measurable subsets of
$\Omega$. We also assume that $\mu$ is absolutely continuous
with respect to Lebesgue measure. If $1\le p< \infty$, let
$E^p_\mu(\Omega)$ denote the class of locally Lebesgue integrable
functions on $\Omega$ with distributional derivatives in
$L^p_\mu(\Omega)$. If $1\le N \le \infty$, we say that a set $Y
\subset L^N_w(\Omega)\cap E^p_\mu(\Omega)$ (intersection of function
spaces instead of normed spaces of equivalence classes) is {\it
bounded in} $L^N_w(\Omega)\cap E^p_\mu(\Omega)$ if
\[
\sup_{f\in Y}\left\{||f||_{L^N_w(\Omega)} + ||\nabla
f||_{L^p_\mu(\Omega)}\right\} < \infty.
\]

We use $D$ to denote a generic open Euclidean ball. The radius and
center of $D$ will be denoted $r(D)$ and $x_D$, and if $C$ is a
positive constant, $CD$ will denote the ball concentric with $D$ whose
radius is $Cr(D)$.

\begin{thm}\label{main2}
Let $\tilde{\Omega}\subset \Omega$ be open sets in $\mathbb{R}^n$. Let
$w$ be a Borel measure on $\Omega$ with $w(\tilde{\Omega})= w(\Omega)
<\infty$ and $\mu$ be a measure on the Lebesgue measurable sets in
$\Omega$ which is absolutely continuous with respect to Lebesgue
measure. Let $1\le p <\infty$, $1<N\le \infty$ and $\mathscr{S}\subset
L^N_w(\Omega)\cap E^p_\mu(\Omega)$, and suppose that for all
$\epsilon>0$, there exists $\delta_\epsilon>0$ such that
\begin{equation}\label{poincare2}
\|f-f_{D,w}\|\down{L^p_w(D)} \le \epsilon \|\nabla f\|\down{L^p_\mu(D)}
\ \mbox{ for all } f\in \mathscr{S}
\end{equation}
and all Euclidean balls $D$ with $r(D)<\delta_\epsilon$ and $2D\subset
\tilde{\Omega}$. Then for any sequence $\{f_k\} \subset \mathscr{S}$ that
is bounded in $L^N_w(\Omega) \cap E_\mu^p(\Omega)$, there is a
subsequence $\{f_{k_i}\}$ and a function $f \in L^N_w(\Omega)$ such that
$\{f_{k_i}\} \to f$ pointwise a.e.-$w$ in $\Omega$ and in
$L^q_w(\Omega)$ norm for $1\le q<N$.
\end{thm}

Before proving Theorem \ref{main2}, we give typical examples of
$\tilde{\Omega}$ and $w$ with $w(\tilde{\Omega})= w(\Omega)
<\infty$. For any two nonempty sets $E_1, E_2 \subset \mathbb{R}^n$, let
\begin{equation}\label{rho}
\rho(E_1,E_2) = \inf \{|x-y|: x\in E_1, y\in E_2\}
\end{equation}
denote the Euclidean distance between $E_1$ and $E_2$. If $x\in
\mathbb{R}^n$ and $E$ is a nonempty set, we will write $\rho(x,E)$
instead of  $\rho(\{x\},E)$.  Let $\tilde{\Omega}$ be an open subset
of $\Omega$.  If $\Omega$ is bounded and
$\Omega\setminus\tilde{\Omega}$ has Lebesgue measure $0$, 
the measure $w$ on $\Omega$ defined by $dw =\rho(x, \mathbb{R}^n
\setminus \tilde{\Omega})^\alpha dx$ clearly has the desired
properties if $\alpha \ge 0$. The range of $\alpha$ can be increased
to $\alpha >-1$ if $\Omega$ is a Lipschitz domain and $\Omega
\setminus \tilde{\Omega}$ is a finite set. Indeed, if $\partial\Omega$
is described in local coordinates $x= (x_1,\dots,x_n)$ by $x_n= 
F(x_1,\dots,x_{n-1})$ with $F$ Lipschitz, then the distance from $x$
to $\partial\Omega$ is equivalent to $|x_n -F(x_1,\dots, x_{n-1})|$,
and consequently the restriction $\alpha >-1$ guarantees that $w$ is
finite near $\partial\Omega$ by using Fubini's theorem; see also
\cite[Remark 3.4(b)]{C1}. If $\Omega$ is bounded and $\Omega\setminus
\tilde{\Omega}$ is finite, but with no restriction on
$\partial\Omega$, the range can clearly be further increased to
$\alpha >-n$ for the measure $\rho(x, \Omega\setminus
{\tilde{\Omega}})^\alpha dx$. Also note that any $w$ without point
masses satisfies $w(\tilde{\Omega}) = w(\Omega)$ if $\tilde{\Omega}$
is obtained by deleting a countable subset of $\Omega$.

\vspace{.3cm}

\noindent{\bf Proof of Theorem \ref{main2}:} We will verify
the hypotheses of Theorem \ref{general}.   Let
\[
\mathscr{X}(\Omega) = \big\{\g = (g_1,\dots,g_n): |\g|=
\big(\sum_{i=1}^n g_i^2\big)^{1/2} \in L^p_\mu(\Omega)\big\}
\]
and $||\g||_{\mathscr{X}(\Omega)} = ||\g||_{L^p_\mu(\Omega)}$. Then
$$
\|\nabla f\|\down{\mathscr{X}(\Omega)}=\|\nabla
f\|\down{L^p_\mu(\Omega)}\quad\mbox{if $f\in E^p_\mu(\Omega)$.}
$$
If $f \in E^p_\mu(\Omega)$, we may identify $f$ with the pair
$(f,\nabla f)$ since the distributional gradient $\nabla f$ is
uniquely determined by $f$ up to a set of Lebesgue measure zero. Then
$L^N_w(\Omega)\cap E^p_\mu(\Omega)$ can be viewed
as a subset of $L^N_w(\Omega)\times \mathscr{X}(\Omega)$.  In Theorem
\ref{general}, choose ${\cal S}$ to be the particular sequence
$\{f_k\} \subset \mathscr{S}$ in the hypothesis of Theorem \ref{main2},
and choose $\Sigma$ to be the Lebesgue measurable subsets of $\Omega$ and
$\Sigma_0$ to be the collection of balls $D \subset \Omega$. Then
hypotheses (A) and ($B_p$) are valid with $C_2=C_1$ for any $C_1$.
Given $\epsilon >0$, since $w(\tilde{\Omega})=w(\Omega)<\infty$, there
is a compact set $K\subset \tilde{\Omega}$ with $w(\Omega\setminus
K)<\epsilon$. Let $0<\delta'_\epsilon <\rho(K,\R^n\setminus
\tilde{\Omega})$ (where $\rho(K, \R^n\setminus\tilde{\Omega})$ is
interpreted as $\infty$ if $\tilde{\Omega} =\mathbb{R}^n$), let
$\delta_\epsilon$ be as in (\ref{poincare2}), and fix 
$r_\epsilon$  with $0<r_\epsilon < \min\{\delta_\epsilon,
\delta'_\epsilon\}$. By considering the triples of balls in a
maximal collection of pairwise disjoint balls of radius $r_\epsilon/6$
centered in $K$, we obtain a collection
$\{E_{\ell}^\epsilon\}_{\ell}$ of balls of radius $r_\epsilon/2$ which
satisfy $2E_\ell^\epsilon \subset \tilde{\Omega}$, have bounded overlaps
with overlap constant independent of $\epsilon$, and whose union
covers $K$. Since $K$ is compact, we may assume the collection is
finite. Also,
\[ w\big(\Omega \setminus \cup_\ell E_\ell^\epsilon\big) \le
w\big(\Omega \setminus K\big) <\epsilon,
\]
and (\ref{pc}) holds with $F_\ell =E_\ell= E_\ell^\epsilon$ by
(\ref{poincare2}). Theorem \ref{main2} now follows from Theorem
\ref{general} applied to $\Omega$. \qed

In particular, we obtain the following result when $w=\mu$ is a
Muckenhoupt $A_p(\mathbb{R}^n)$ weight, i.e., when $d\mu = dw = \eta
\,dx$ where $\eta(x)$ satisfies
$$
\left(\frac{1}{|D|}\int_D\eta\,dx\right)\left(\frac{1}{|D|}
\int_D \eta^{-1/(p-1)} dx\right)^{p-1} \le C
$$
if $1<p<\infty$ and $|D|^{-1}\int_D \eta\,dx \le C\, \text{essinf}_D
w$ if $p=1$ for all Euclidean balls $D$, with $C$ independent of $D$.
As is well known, such a weight also satisfies the classical doubling
condition
\bea\label{classicaldoubling}
w(D_r(x)) \le C \left(\frac{r}{r'}\right)^\theta w(D_{r'}(x)), \quad
0<r'<r<\infty,
\eea
with $\theta \ge np-\epsilon$ for some $\epsilon>0$ if $p>1$, and with
$\theta =n$ if $p=1$, where $C$ and $\theta$ are independent of $r,
r', x$. 

We denote by $W^{1,p,w}(\Omega)$ the weighted Sobolev space defined as
all functions in $L^p_w(\Omega)$ whose distributional gradient is in
$L^p_w(\Omega)$. Thus $W^{1,p,w}(\Omega) = L^p_w(\Omega) \cap
E^p_w(\Omega)$.  If $w(\Omega)<\infty$, it follows that $L^N_w(\Omega)\cap
E^p_w(\Omega) \subset W^{1,p,w}(\Omega)$ when $N\ge p$, and that the
opposite containment holds when $N\le p$.

\begin{thm}\label{nondegen} Let $1\le p <\infty$, $w\in
A_p(\mathbb{R}^n)$ and $\Omega$ be an open set in $\mathbb{R}^n$ with
$w(\Omega)<\infty$.  If $1< N \le \infty$, then any bounded subset of
$L^N_w(\Omega)\cap E^p_w(\Omega)$ is precompact in $L^q_w(\Omega)$ if
$1\le q <N$. Consequently, if $N>p$ and $\mathscr{S}$ is a subset of
$W^{1,p,w}(\Omega)$ with
\begin{equation}\label{wsob}
\|f\|_{L^N_w(\Omega)}\le C(\|f\|_{L^p_w(\Omega)}+\|\nabla
f\|_{L^p_w(\Omega)}) \ \mbox{ for all } f\in \mathscr{S},
\end{equation}
then any set in $\mathscr{S}$ that is bounded in $W^{1,p,w}(\Omega)$
is precompact in $L^q_w(\Omega)$ for $1\le q<N$.

If $\Omega$ is a John domain, there exists $N>p$ ($N$ can be $\theta
p/(\theta -p)$ for some $\theta>p$ as described after
(\ref{classicaldoubling})) such that $W^{1,p,w}(\Omega)$ is compactly
embedded in $L^q_w(\Omega)$ for  $1\le q <N$. In particular, the
embedding of $W^{1,p,w}(\Omega)$ into $L^p_w(\Omega)$ is compact when
$w\in A_p(\mathbb{R}^n)$ and $\Omega$ is a John domain.

\end{thm}

\begin{remark}

When $w=1$ and $p<n$, the choices $N= np/(n-p)$ and $\mathscr{S} =
W_0^{1,p}(\Omega)$  guarantee (\ref{wsob}) by the classical Sobolev
inequality for functions in $W_0^{1,p}(\Omega)$ (see
e.g. \cite[Theorem 7.10]{GT}); here $W_0^{1,p}(\Omega)$ denotes the
closure in $W^{1,p}(\Omega)$ of the class of Lipschitz functions with
compact support in $\Omega$. Consequently, the classical
Rellich-Kondrachov theorem giving the compact embedding of
$W_0^{1,p}(\Omega)$ in $L^q(\Omega)$ for $1\le q <np/(n-p)$ follows as
a special case of the first part of Theorem \ref{nondegen}.
\end{remark}
\medskip

{\bf Proof. }  We will apply Theorem \ref{main2} with $w=\mu$. Fix $p$
and $w$ with $1\le p <\infty$ and $w\in A_p(\mathbb{R}^n)$. By
\cite{FKS}, there is a constant $C$ such that the weighted Poincar\'e
inequality
\[
||f -f_{D,w}||_{L^p_w(D)} \le C r(D) ||\nabla f||_{L^p_w(D)}, \quad f\in
  C^\infty(\Omega),
\]
holds for all Euclidean balls $D\subset \Omega$. Then since
$C^\infty(\Omega)$ is dense in $L^N_w(\Omega)\cap E^p_w(\Omega)$ if
$1\le N <\infty$ (see e.g. \cite{Tur}), by fixing any $\epsilon>0$ we
obtain from Fatou's lemma that for all balls $D \subset \Omega$ with
$C r(D) \le \epsilon$,
\[
||f -f_{D,w}||_{L^p_w(D)} \le \epsilon\, ||\nabla f||_{L^p_w(D)}
  \quad\text{if $f\in L^N_w(\Omega)\cap E^p_w(\Omega)$}.
\]
The same holds when $N=\infty$ since $L^\infty_w(\Omega) =
L^\infty(\Omega) \subset L^p_w(\Omega)$ due to the assumptions $w\in
A_p(\mathbb{R}^n)$ and $w(\Omega)<\infty$. With $1<N\le \infty$, the
first statement of the theorem now follows from Theorem \ref{main2},
and the second statement is a corollary of the first one.

Next, let $\Omega$ be a John domain. Choose $\theta>p$ so that $w$
satisfies (\ref{classicaldoubling}) and define $N =
\theta p/(\theta-p)$. Then $N>p$ and by \cite[Theorem 1.8 (b) or
Theorem 4.1]{CW1},
\[
||f -f_{\Omega,w}||_{L^N_w(\Omega)} \le C ||\nabla
  f||_{L^p_w(\Omega)}, \quad \forall f\in C^\infty(\Omega).
\]
Again, the inequality remains true for functions in
$W^{1,p,w}(\Omega)$ by density and Fatou's lemma. It is now clear that
(\ref{wsob}) holds, and the last part of the theorem follows.\qed
\\

Our next example involves domains in $\mathbb{R}^n$ which are more
restricted. For special $\Omega$, there are values $N>1$ such that
\begin{equation}\label{embedding}
\|f\|\down{L^N(\Omega)}\le C\big(\|f\|\down{L^1(\Omega)}+\|\nabla
f\|\down{L^p(\Omega)}\big)
\end{equation}
for all $f\in L^1(\Omega)\cap E^p(\Omega)$. Note that
if $\Omega$ has finite Lebesgue measure, then $W^{1,p}(\Omega) \subset
L^1(\Omega)\cap E^p(\Omega)$. As we will explain, (\ref{embedding}) is
true for some $N>1$ if $\Omega$ is an $s$-John domain in $\mathbb{R}^n$
and $1 \le s <1+\frac{p}{n-1}$. Recall that for $1\le s<\infty$, a
bounded domain $\Omega \subset \mathbb{R}^n$ is called an $s$-John
domain with central point $x'\in \Omega$ if for some constant $c>0$
and all $x \in \Omega$ with $x \neq x'$, there is a curve $\Gamma:
[0,l] \rightarrow \Omega$ so that $\Gamma(0)=x, \Gamma(l) =x'$,
\[
|\Gamma(t_1)-\Gamma(t_2)| \le t_2-t_1\quad\mbox{for all $[t_1,t_2]
 \subset [0,l]$, and}
\]
\[
\rho(\Gamma(t), \Omega^c) \ge c\, t^s\quad\mbox{for all $t\in [0,l]$.}
\]
The terms $1$-John domain and John domain are the same. When $\Omega$
is an $s$-John domain for some $s\in [1, 1+p/(n-1))$, it is shown in
\cite{KM}, \cite{CW1}, \cite{CW2} that (\ref{embedding}) holds for all
finite $N$ with
\begin{equation}\label{scondition}
\frac{1}{N} \ge \frac{s(n-1)-p+1}{np}
\end{equation}
and for all $f\in W^{1,p}(\Omega)$ without any support restrictions.
Note that the right side of (\ref{scondition}) is strictly less than
$1/p$ for such $s$, and consequently there are values $N>p$ which
satisfy (\ref{scondition}). For $N$ as in (\ref{scondition}), the
global estimate
\begin{equation}\label{KMpoincare}
||f-f_\Omega||_{L^{N}(\Omega)} \le C ||\nabla f||_{L^p(\Omega)},
  \quad f_\Omega = \int_\Omega f(x) dx/|\Omega|,
\end{equation}
is shown to hold if $f\in Lip_{loc}(\Omega)$ in \cite{CW2}, and then
follows for all $f\in L^1(\Omega)\cap E^{p}(\Omega)$; see the proof of
Theorem \ref{weightedsjohn} for related comments.  Inequality
(\ref{embedding}) is clearly a consequence of (\ref{KMpoincare}).

More generally, weighted versions of (\ref{KMpoincare}) hold for
$s$-John domains and lead to weighted compactness results, as we
now show.  Let $1\le p< \infty$, and for real $\alpha$ and
$\rho(x,\Omega^c)$ as in (\ref{rho}), let $L^p_{\rho^\alpha
dx}(\Omega)$ be the class of Lebesgue measurable $f$ on $\Omega$ with
\[
||f||_{L^p_{\rho^\alpha dx}(\Omega)} = \left(\int_\Omega |f(x)|^p
  \rho(x,\Omega^c)^\alpha dx\right)^{1/p} <\infty.
\]

\begin{thm}\label{weightedsjohn}  Suppose that $1\le s <\infty$
and $\Omega$ is an $s$-John domain in $\mathbb{R}^n$. Let $p, a, b$
satisfy $1\le p<\infty$, $a\ge 0$, $b\in \mathbb{R}$ and $b-a<p$.

(i) If
\begin{equation}\label{scondition1}
n+a>s(n-1+b)-p+1,
\end{equation}
then for any $1\le q <\infty$ such that
\begin{equation}\label{scondition2}
\frac{1}{q} > \max\left\{\frac{1}{p}-\frac{1}{n},
\frac{s(n-1+b)-p+1}{(n+a)p}\right\},
\end{equation}
$L^1_{\rho^a dx}(\Omega)\cap E^p_{\rho^b dx}(\Omega)$ is compactly
embedded in $L^q_{\rho^a dx}(\Omega)$.

(ii) If $p>1$ and
\begin{equation}\label{s1condition}
n+ap>s(n-1+b)-p+1\ge n+a,
\end{equation}
then for any $1\le q<\infty$ such that
\begin{equation}\label{s2condition}
\frac{a}{q} > \max\left\{\frac{b}{p}-1,
\frac{s(n-1+b)-p-n+1}{p}\right\},
\end{equation}
$L^1_{\rho^a dx}(\Omega)\cap E^p_{\rho^b dx}(\Omega)$ is compactly
embedded in $L^q_{\rho^a dx}(\Omega)$.
\end{thm}

\begin{remark}\label{re-sjohn}
\begin{enumerate}
\item If $a=b=0$, (\ref{scondition1}) is the same as $s<
1+\frac{p}{n-1}$. If $a=0$, (\ref{s1condition}) never holds.
\item The requirement that $b-a<p$ follows from (\ref{scondition1})
and (\ref{scondition2}) by considering the cases $n-1+b \ge 0$ and
$n-1+b <0$ separately. Hence $b-a<p$ automaticallly holds in part (i),
but it is an assumption in part (ii). Also, (\ref{s1condition}) and
(\ref{s2condition}) imply that $q<p$, and consequently that $p>1$.
\item Conditions (\ref{scondition1}) and (\ref{scondition2}) imply
there exists $N \in (p, \infty)$ with
\begin{equation}\label{scondition3}
\frac{1}{q}>\frac{1}{N} > \max\left\{\frac{1}{p}-\frac{1}{n},
\frac{s(n-1+b)-p+1}{(n+a)p}\right\}.
\end{equation}
Conversely, (\ref{scondition1}) holds if there exists $N\in (p,
\infty)$ so that (\ref{scondition3}) holds.
\item Assumption (\ref{s2condition}) ensures that there exists $N\in
(q, \infty)$ such that (\ref{s2condition}) holds with $q$ replaced by
$N$.
\end{enumerate}

\end{remark}

\noindent{\bf Proof:} This result is also a consequence of Theorem
\ref{main2}, but we will deduce it from Theorem \ref{general} by
using arguments like those in the proofs of Theorems \ref{main2}
and \ref{nondegen}. Fix $a,b,p,q$ as in the hypothesis and denote
$\rho(x) = \rho(x,\Omega^c)$. Choose $w= \rho^a dx$ and note that
$w(\Omega) <\infty$ since $a\ge 0$ and $\Omega$ is now bounded. Define
\[
\mathscr{X}(\Omega) = \big\{\g = (g_1,\dots,g_n): |\g|\in
L^p_{\rho^b dx}(\Omega)\big\}
\]
and $||\g||_{\mathscr{X}(\Omega)} = ||\g||_{L^p_{\rho^a
dx}(\Omega)}$. Fix $\epsilon>0$ and choose a compact set $K \subset
\Omega$ with $|\Omega\setminus K|_{\rho^a dx} :=
\int_{\Omega\setminus K} \rho^a dx < \epsilon$. Also choose
$\delta_\epsilon'$ with $0<\delta_\epsilon' < \rho(K, \Omega^c)$,
where $\rho(K, \Omega^c)$ is the Euclidean distance between $K$ and
$\Omega^c$.

If $D$ is a Euclidean ball with center $x_D\in K$ and
$r(D) < \frac{1}{2}\delta_\epsilon'$, then
$2D\subset \Omega$ and $\rho(x)$ is essentially constant on $D$; in
fact, for such $D$,
\[
\frac{1}{2}\rho(x_D) \le \rho(x) \le \frac{3}{2}\rho(x_D),\quad x \in
D.
\]

We claim that for such $D$, the simple unweighted Poincar\'e estimate
\[
||f -f_D||_{L^p(D)} \le C r(D) ||\nabla f||_{L^p(D)}, \quad f\in
  Lip_{loc}(\Omega),
\]
where $f_D= f_{D,dx}$, implies that for $f\in Lip_{loc}(\Omega)$,
\begin{equation}\label{weightedpc}
||f -f_{D,\rho^a dx}||_{L^p_{\rho^a dx}(D)} \le {\tilde C}
  \big(r(D)^{\frac{a-b}{p}}  + \diam(\Omega)^{\frac{a-b}{p}}\big) r(D)
  ||\nabla f||_{L^p_{\rho^b dx}(D)},
\end{equation}
where $f_{D,\rho^a dx}= \int_D f\rho^a dx/\int_D \rho^a dx$ and
${\tilde C}$ depends on $C,a,b$ but is independent of $D, f$. To show
this, first note that for such $D$, since $\rho\sim \rho(x_D)$ on $D$,
the simple Poincar\'e estimate immediately gives
\[
||f -f_D||_{L^p_{\rho^a dx}(D)} \le {\tilde C} \rho(x_D)^{\frac{a
-b}{p}} r(D) ||\nabla f||_{L^p_{\rho^b dx}(D)}, \quad f\in
Lip_{loc}(\Omega),
\]
and then a similar estimate with $f_D$ replaced by $f_{D,\rho^a dx}$
follows by standard arguments. Clearly (\ref{weightedpc}) will now
follow if we show that
\[
\rho(x_D)^{\frac{a-b}{p}}\le r(D)^{\frac{a-b}{p}}
  + \diam(\Omega)^{\frac{a-b}{p}}\quad \text{for such $D$.}
\]
However, this is clear since $r(D)\le \rho(x_D) \le diam(\Omega)$ for
$D$ as above, and (\ref{weightedpc}) is proved.

We can now apply the weighted density result of \cite{H}, \cite{HK1} to
conclude that (\ref{weightedpc}) holds for all $f \in
L^1_{\rho^a dx}(\Omega)\cap E^p_{\rho^bdx}(\Omega)$ and all balls $D$
with $x_D\in K$ and $r(D) <\frac{1}{2}\delta_\epsilon'$.

Recall that $\frac{a-b}{p} +1>0$. Thus there exists $r_\epsilon$ with
$0 <r_\epsilon <\frac{1}{2}\delta_\epsilon'$ and
\[
{\tilde C} \big(r_\epsilon^{\frac{a-b}{p}}
  + \diam(\Omega)^{\frac{a-b}{p}}\big) r_\epsilon  < \epsilon.
\]
Let $\Sigma$ and $\Sigma_0$ be as in the proof of Theorem \ref{main2},
and let $\{E_\ell\}_\ell=\{F_\ell\}_\ell$ be the triples of balls in a
maximal collection of pairwise disjoint balls centered in $K$ with
radius $\frac{1}{3}r_\epsilon$. Then (\ref{weightedpc}) and the choice
of $r_\epsilon$ give the desired version of (\ref{pc}), namely
\[
||f -f_{D,\rho^a dx}||_{L^p_{\rho^a dx}(D)} \le \epsilon ||\nabla
  f||_{L^p_{\rho^b dx}(D)}
\]
for $D=E_\ell$ and $f\in L^1_{\rho^a dx}(\Omega)\cap
E^p_{\rho^bdx}(\Omega)$. Next, use the last two parts of Remark
\ref{re-sjohn} to choose $N\in (q,\infty)$ so that either
(\ref{scondition2}) or (\ref{s2condition}) holds with $q$ there
replaced by $N$. Every $f\in L^1_{\rho^a dx}(\Omega)\cap E^p_{\rho^b
dx}(\Omega)$ then satisfies the global Poincar\'e estimate
\begin{equation} \label{globalsjohn}
||f- f_{\Omega, \rho^a dx}||_{L^N_{\rho^a dx}(\Omega)} \le C ||\nabla
  f||_{L^p_{\rho^b dx}(\Omega)}, \quad f\in L^1_{\rho^a
  dx}(\Omega)\cap E^p_{\rho^b dx}(\Omega),
\end{equation}
where $f_{\Omega,\rho^a  dx} = \int_\Omega f \,\rho^a dx/\int_\Omega
\rho^a dx$. In fact, under the hypothesis of Theorem
\ref{weightedsjohn}, this is proved for $f\in Lip_{loc}(\Omega) \cap
L^1_{\rho^a dx}(\Omega)\cap E^p_{\rho^bdx}(\Omega)$ in \cite {CW2} for
example, and then follows for all $f\in L^1_{\rho^a dx}(\Omega)\cap
E^p_{\rho^bdx}(\Omega)$ by the density result of \cite{H}, \cite{HK1}
and Fatou's lemma. By (\ref{globalsjohn}),
\[
||f ||_{L^N_{\rho^a dx}(\Omega)} \le C ||f||_{L^1_{\rho^a
  dx}(\Omega)} + C ||\nabla f||_{L^p_{\rho^b dx}(\Omega)}
\]
for the same class of $f$. The remaining details of the proof are
left to the reader. \qed

In passing, we mention that the role played by the distance function
$\rho(x,\Omega^c)$ in Theorem \ref{weightedsjohn} can instead be
played by
$$
\rho_0(x) = \inf\{|x-y|: y\in \Omega_0\},\quad x\in \Omega,
$$
for certain $\Omega_0 \subset \Omega^c$; see \cite[Theorem 1.6]{CW2}
for a description of such $\Omega_0$ and the required Poincar\'e
estimate, and note that the density result in \cite{HK1} holds for
positive continuous weights.

\section{Applications in the Degenerate Case}
\setcounter{equation}{0}

\setcounter{theorem}{0}

In this section, $\Omega$ denotes a fixed open set in $\mathbb{R}^n$,
possibly unbounded. For $(x,\xi) \in \Omega\times \mathbb{R}^n$, we
consider a nonnegative quadratic form $\xi' Q(x)\xi$ which may
degenerate, i.e., which may vanish for some $\xi \neq 0$. Such
quadratic forms occur naturally in the context of subelliptic
equations and give rise to degenerate Sobolev spaces as
discussed below. Our goal is to apply Theorem \ref{general} to obtain
compact embedding of these degenerate spaces into Lebesgue spaces
related to the gain in integrability provided by Poincar\'e-Sobolev
inequalities. The framework that we will use contains the subelliptic
one developed in \cite[2]{SW1}, where regularity theory for weak
solutions of linear subelliptic equations of second order in
divergence form is studied.

\subsection{Standing Assumptions}

We now list some notation and assumptions that will be in force
everywhere in \S 3 even when not explicitly mentioned.

\begin{defn} A function $d$ is called a finite symmetric
quasimetric (or simply a quasimetric) on $\Omega$ if
$d:\Omega\times\Omega\ra [0,\infty)$ and there is a constant
$\kappa\geq 1$ such that for all $x,y,z\in\Omega$,
\bea
d(x,y) &=& d(y,x),\nonumber\\
d(x,y) &=& 0 \iff x=y, \textrm{ and} \nonumber\\
\label{tri} d(x,y)&\leq&\kappa [d(x,z)+d(z,y)].
\eea
\end{defn}
If $d$ is a quasimetric on $\Omega$, we refer to the pair
$(\Omega,d)$ as a quasimetric space.  In some applications, $d$ is
closely related to $Q(x)$.  For example, $d$ is sometimes chosen to be
the Carnot-Carath\'eodory control metric related to $Q$; cf. \cite{SW1}.

Given $x\in\Omega$, $r>0$, and a quasimetric $d$, the subset of
$\Omega$ defined by
\bea B_r(x) = \{y\in\Omega\;:\; d(x,y) <r\}\nonumber
\eea
will be called the quasimetric $d$-ball centered at $x$ of radius
$r$. Note that every $d$-ball $B = B_r(x)$ satisfies $B \subset
\Omega$ by definition.

It is sometimes possible, and desirable in case the boundary of
$\Omega$ is rough, to be able to work only with $d$-balls that are
deep inside $\Omega$ in the sense that their Euclidean closures
$\overline{B}$ lie in $\Omega$. See part (ii) of Remark \ref{various}
for comments about being able to use such balls.

Recall that $D_s(x)$ denotes the ordinary Euclidean ball of radius $s$
centered at $x$. We always assume that $d$ is related as follows to
the standard Euclidean metric:
\bea\label{new2}
\mbox{$\forall\, x\in\Omega$ and $r>0$, $\exists \, s=s(x,r)>0$ so
that } D_{s}(x) \subset B_r(x).
\eea

\begin{remark}
Condition (\ref{new2}) is clearly true if $d$-balls are open, and it
is weaker than the well-known condition of C. Fefferman and Phong
stating that for each compact $K\subset\Omega$, there are constants
$\beta, r_0>0$ such that $D_{r^\beta}(x)\subset B_r(x)$ for all $x \in
K$ and $0<r<r_0$.
\end{remark}

Throughout \S 3, $Q(x)$ denotes a fixed Lebesgue measurable $n\times
n$ nonnegative symmetric matrix on $\Omega$ and we assume that every
$d$-ball $B$ centered in $\Omega$ is Lebesgue measurable.  We will deal
with three locally finite measures $w,\nu,\mu$ on the Lebesgue
measurable subsets of $\Omega$, each with a particular role.  In \S
3.3, where only global results are developed, we will assume
$w(\Omega)<\infty$ but this assumption is not required for the local
results of \S 3.4.  The measure $\mu$ is assumed to be absolutely
continuous with respect to Lebesgue measure; the comment following
(\ref{sobnorm}) explains why this assumption is natural. In \S 3, we
sometimes assume that $w$ is absolutely continuous with respect to $\nu$,
but we drop this assumption completely in the Appendix.

We do not require the existence of a doubling measure for the
collection of $d$-balls, but we always assume that $(\Omega,d)$
satisfies the weaker local geometric doubling property given in the
next definition; see \cite{HyM} for a global version.

\begin{defn}\label{doublingdef} A quasimetric space $(\Omega,d)$
satisfies the {\it local geometric doubling condition} if for every
compact $K\subset\Omega$, there exists $\delta'=\delta'(K)>0$ such
that for all $x\in K$ and all $0<r'<r < \delta'$, the number of
disjoint $d$-balls of radius $r'$ contained in $B_{r}(x)$ is at most
a constant ${\cal C}_{r/r'}$ depending on $r/r'$ but not on $K$.
\end{defn}

\subsection{Degenerate Sobolev Spaces $\wa{p}{\Omega},\;\wb{p}{\Omega}$}

We will define weighted degenerate Sobolev spaces by using an approach
like the one in \cite{SW2} for the unweighted case. We first define an
appropriate space of vectors, including vectors which will eventually
play the role of gradients, where size is measured relative to the
nonnegative quadratic form
 \bea Q(x,\xi) = \xi'Q(x)\xi,\quad (x,\xi)\in
\Omega\times\mathbb{R}^n. \nonumber
\eea
For $1\leq p<\infty$, consider the collection of measurable
$\mathbb{R}^n$-valued functions $\vec{g}(x) =  (g_1(x),...,g_n(x))$
satisfying
\bea \label{ellpnorm}||\vec{g}||_{{\cal L}^p_\mu(\Omega,Q)} =
\Big\{\int_\Omega Q(x,\vec{g}(x))^\frac{p}{2}d\mu\Big\}^\frac{1}{p} =
 \Big\{\int_\Omega |\sqrt{Q(x)}\vec{g}(x)|^pd\mu\Big\}^\frac{1}{p}
 <\infty.
\eea
We identify any two functions $\vec{g}, \vec{h}$ in the collection for
which $||\vec{g} - \vec{h}||_{{\cal L}^p_\mu(\Omega,Q)} = 0$. Then
(\ref{ellpnorm}) defines a norm on the resulting space of equivalence
classes. The form-weighted space ${\cal  L}^p_\mu(\Omega,Q)$ is
defined to be the collection of these equivalence classes, with norm
(\ref{ellpnorm}). By using methods similar to those in \cite{SW2}, it
follows that ${\cal L}^2_\mu(\Omega,Q)$ is a Hilbert space and ${\cal
L}^p_\mu(\Omega,Q)$ is a Banach space for $1\le p<\infty$.

Now consider the (possibly infinite) norm on $Lip_{loc}(\Omega)$
defined by
\bea \label{sobnorm}||f||_{\wa{p}{\Omega}} =
||f||_{L^p_\nu(\Omega)} + ||\nabla f||_{{\cal L}^p_\mu(\Omega,Q)}.
\eea
We comment here that our standing assumption that $\mu(Z)=0$ when $Z$
has Lebesgue measure $0$ assures that $||\nabla f||_{{\cal
L}^p_\mu(\Omega,Q)}$ is well-defined if $f\in Lip_{loc}(\Omega)$; in
fact, for such $f$, the Rademacher-Stepanov theorem implies that
$\nabla f$ exists a.e. in $\Omega$ with respect to Lebesgue measure.

\begin{defn}\label{spacedef} Let $1\leq p <\infty$.
\begin{enumerate}
\item The degenerate Sobolev space $\wa{p}{\Omega}$ is the completion
under the norm (\ref{sobnorm}) of the set
$$
Lip_{Q,p}(\Omega) = Lip_{Q,p,\nu,\mu}(\Omega) = \{ f\in
Lip_{loc}(\Omega)\;:\; ||f||_{\wa{p}{\Omega}}<\infty\}.
$$
\item The degenerate Sobolev space $\wb{p}{\Omega}$ is the completion
under the norm (\ref{sobnorm}) of the set $Lip_{Q,p,0}(\Omega) =
Lip_0(\Omega)\cap Lip_{Q,p}(\Omega)$, where $Lip_0(\Omega)$ denotes the
collection of Lipschitz functions with compact support in $\Omega$. If
$Q \in L^{p/2}_{loc}(\Omega)$, then $Lip_{Q,p,0}(\Omega) =
Lip_0(\Omega)$ since $\nu$ and $\mu$ are locally finite.
\end{enumerate}
\end{defn}

We now make some comments about $\wa{p}{\Omega}$, most
of which have analogues for $\wb{p}{\Omega}$. By definition,
$\wa{p}{\Omega}$ is the Banach space of equivalence classes of
Cauchy sequences of $Lip_{Q,p}(\Omega)$ functions with respect to
the norm (\ref{sobnorm}).  Given a Cauchy sequence $\{f_j\}$ of
$Lip_{Q,p}(\Omega)$ functions, we denote its equivalence class by
$[\{f_j\}]$.  If $\{v_j\}\in [\{f_j\}]$, then $\{v_j\}$ is a Cauchy
sequence in $L^p_\nu(\Omega)$ and $\{\nabla v_j\}$ is a Cauchy
sequence in ${\cal L}^p_\mu(\Omega,Q)$.  Hence, there is a pair
$(f,\vec{g}) \in L^p_\nu(\Omega) \times {\cal L}^p_\mu(\Omega,Q)$ so
that
\bea
||v_j-f||_{L^p_\nu(\Omega)} \ra 0\quad \textrm{and }||\nabla v_j -
\vec{g}||_{{\cal L}^p_\mu(\Omega,Q)}\ra 0\nonumber
\eea
as $j\ra \infty$.  The pair $(f,\vec{g})$ is uniquely determined by
the equivalence class $[\{f_j\}]$, i.e., is independent of a
particular $\{v_j\}\in [\{f_j\}]$. We will say that $(f,\vec{g})$ is
{\it represented by} $\{v_j\}$. We obtain a Banach space isomorphism ${\cal
J}$ from $\wa{p}{\Omega}$ onto a closed subspace $\wc{p}{\Omega}$ of
$L^p_\nu(\Omega)\times {\cal L}^p_\mu(\Omega,Q)$ by setting \bea
{\cal J}([\{f_j\}]) = (f,\vec{g}).
\eea
We will often not distinguish between $\wa{p}{\Omega}$ and
$\wc{p}{\Omega}$. Similarly, $\we{p}{\Omega}$ will denote the image of
$\wb{p}{\Omega}$ under ${\cal J}$, but we often consider these spaces
to be the same.

It is important to think of a typical element of $\wc{p}{\Omega}$, or
$\wa{p}{\Omega}$, as a pair $(f,\vec{g})$ as above, and not simply as
the first component $f$. In fact, if $(f,\vec{g}) \in \wc{p}{\Omega}$,
the vector $\vec{g}$ may not be uniquely determined by $f$; see
\cite[Section 2.1]{FKS} for a well known example.

If $f\in Lip_{Q,p}(\Omega)$, then the pair $(f,\nabla f)$ may be
viewed as an element of $W_{\nu,\mu}^{1,p}(\Omega,Q)$ by identifying it
with the equivalence class $[\{f\}]$ corresponding to the sequence
each of whose entries is $f$.  When viewed as a class, $(f, \nabla f)$
generally contains pairs whose first components are not Lipschitz
functions; for example, if $f \in Lip_{Q,p}(\Omega)$ and $F$ is any
function with $F=f$ a.e.-$\nu$, then $(f, \nabla f) = (F, \nabla f)$
in $W^{1,p}_{\nu,\mu}(\Omega,Q)$. However, in what follows, when we
consider a pair $(f, \nabla f)$ with $f \in Lip_{Q,p}(\Omega)$,
we will {\it not} adopt this point of view. Instead we will identify
an $f\in Lip_{Q,p}(\Omega)$ with the single pair $(f, \nabla f)$ whose
first component is $f$ (defined everywhere in $\Omega$) and whose
second component is $\nabla f$, which exists a.e. with respect to Lebesgue
measure by the Rademacher-Stepanov theorem. This convention lets us
avoid assuming that $w$ is absolutely continuous with respect to
$\nu$, written $w<<\nu$, in Poincar\'e-Sobolev estimates for
$Lip_{Q,p}(\Omega)$ functions. We will reserve the notation ${\cal
H}$ for subsets of $Lip_{Q,p}(\Omega)$ viewed in this way.

On the other hand, ${\cal W}$ will denote various subsets of
$W^{1,p}_{\nu,\mu}(\Omega,Q)$ with elements viewed as equivalence
classes. When our hypotheses are phrased in terms of such ${\cal W}$,
we will assume that $w<<\nu$ in order to avoid technical difficulty
associated with sets of measure $0$; see the comment after
(\ref{poincare*}). In the Appendix, we drop the assumption $w<<\nu$
altogether.

We will abuse the notation (\ref{sobnorm}) by writing
\bea\label{pairnormlip}
||(f,\nabla f)||_{\wa{p}{\Omega}} = ||f||_{L^p_\nu(\Omega)} + ||\nabla
  f||_{{\cal L}^p_\mu(\Omega,Q)},\quad f \in   Lip_{Q,p}(\Omega),
\eea
and we extend this to generic $(f,\vec{g})\in W^{1,p}_{\nu,
\mu}(\Omega, Q)$ by writing
\bea\label{pairnorm}
||(f,\vec g)||_{\wa{p}{\Omega}} = ||f||_{L^p_\nu(\Omega)} +
 ||\vec{g}||_{{\cal L}^p_\mu(\Omega,Q)}.
\eea

\subsection{Global Compactness Results for Degenerate Spaces}

In this section, we state and prove compactness results which apply to
the entire set $\Omega$. Results which are more local are given in \S 3.4.

In order to apply Theorem \ref{general} in this setting, we will use
the following version of Poincar\'e's inequality for $d$-balls.

\begin{defn}\label{poincaredef} Let $1\leq p<\infty$,
$Lip_{Q,p}(\Omega)$ be is as in Definition \ref{spacedef}, and ${\cal
H} \subset Lip_{Q,p}(\Omega)$.  We say that the \emph{Poincar\'e
property of order $p$ holds for} ${\cal H}$ if there is a constant
$c_0\ge 1$ so that for every $\epsilon>0$ and every compact set
$K\subset\Omega$, there exists $\delta = \delta(\epsilon,K)>0$ such
that for all $f \in {\cal H}$ and every $d$-ball $B_r(y)$ with $y\in
K$ and $0<r< \delta$, 
\bea \label{poincare} \left(\int_{B_r(y)} |f-f_{B_r(y),w}|^{p} dw
\right)^\frac{1}{p} \leq
\epsilon ||(f,\nabla f)||_{W_{\nu,\mu}^{1,p}(B_{c_0r}(y),Q)}.
\eea
\end{defn}

\begin{remark} \label{various}

(i) Inequality (\ref{poincare}) is not of standard Poincar\'e
form.  A more typical form is
\bea \label{poincare3}
\left(\fracc{1}{w(B_r(y))}\int_{B_r(y)} |f-f_{B_r(y),w}|^{p}
dw\right)^\frac{1}{p}\hspace{2in} \nonumber\\
\leq Cr \left(\fracc{1}{\mu(B_{c_0r}(y))}\int_{B_{c_0r}(y)}
|\sqrt{Q}\nabla f|^{p}d\mu\right)^\frac{1}{{p}}.
\eea
In \cite[2]{SW1} and \cite{R1}, the unweighted version of
\eqref{poincare3} with $p=2$ is used. Let $\rho(x,\partial\Omega)$ and
$\rho(E,\partial\Omega)$ be as in (\ref{rho}). In \cite{SW2}, the
unweighted form of (\ref{poincare3}) with $p=2$ is assumed for all $f\in
Lip_{Q,2}(\Omega)$ and all $B_r(y)$ with $y\in \Omega$ and $0<r<
\delta_0 \rho(y,\partial\Omega)$ for some $\delta_0\in (0,1)$
independent of $y,r$. If $K$ is a compact set in $\Omega$, this
version would then hold for all $B_r(y)$ with $y\in K$ and $0<r<
\delta_0 \rho(K,\partial\Omega)$. For general $p, w$ and $\mu$, if for
every compact $K\subset\Omega$, (\ref{poincare3}) is valid for all
$B_r(y)$ with $y\in K$ and $0<r< \delta_0 \rho(K, \partial\Omega)$, then
(\ref{poincare}) follows easily provided
\bea \label{balancing}
\lim_{r\ra 0} \left\{\sup_{y\in K}r^{p}
\fracc{w(B_r(y))}{\mu(B_{c_0r}(y))}\right\} = 0
\eea
for every compact $K\subset\Omega$. Note that (\ref{balancing})
automatically holds if $w=\mu$.

If both (\ref{poincare3}) and (\ref{balancing}) hold, then
(\ref{poincare}) is true for any choice of $\nu$. In this situation,
one can pick $\nu =w$ in order to avoid technicalities encountered
below when $w$ is not absolutely continuous with respect to
$\nu$.

(ii) Especially when $\partial\Omega$ is rough, it is simplest to deal
only with $d$-balls $B$ which stay away from $\partial\Omega$, i.e.,
which satisfy
\bea \label{closure}
\overline{B} \subset \Omega.
\eea
We can always assume this for the balls in (\ref{poincare}) if
the converse of (\ref{new2}) is also true, namely if
\bea \label{eucled2}
\forall \, x\in \Omega \ \mbox{ and } r>0, \ \exists \,s= s(r,x)>0 \
\mbox{ such that } B_s(x)\subset D_r(x).
\eea
To see why, let us first show that given a compact set $K$ and an open
set $G$ with $K\subset G\subset\Omega$, there exists $t>0$ so
that $\overline{B_t(y)} \subset G$ for all $y\in K$. Indeed, for such
$K$ and $G$, let  
$t'=\frac{1}{2}\rho(K,G^c)$. By (\ref{eucled2}), for each $x\in
K$ there exists $r(x)>0$ so that $B_{r(x)}(x)\subset
D_{t'}(x)$. Further, by (\ref{new2}), there exists $s(x)>0$ so that
$D_{s(x)}(x)\subset B_{r(x)/(2\k)}(x)$, where $\kappa$ is as in
(\ref{tri}). Since $K$ is compact, we may choose finite collections
$\{B_{r_i/(2\k)}(x_i)\}$ and $\{D_{s_i}(x_i)\}$ with $x_i \in K$,
$r_i = r(x_i)$, $s_i = s(x_i)$, and $K\subset \bigcup
D_{s_i}(x_i) \subset \bigcup B_{r_i/(2\k)}(x_i)$.  Now set $t=\min\{
r_i/(2\k)\}$. Let $y\in K$ and choose $i$ such that $y\in
B_{r_i/(2\k)}(x_i)$. By (\ref{tri}), $B_t(y)\subset B_{r_i}(x_i)$ and
consequently $B_t(y)\subset D_{t'}(x_i)$. Since 
$\overline{D_{t'}(x_i)} \subset G$, we obtain $\overline{B_{t}(y)}
\subset G$ for every $y\in K$, as desired. In particular,
$\overline{B_t(y)} \subset \Omega$ for all $y\in K$. Since the validity of
(\ref{poincare}) for some $\delta=\delta(\epsilon, K)$ implies its
validity for min\,$\{\delta,t\}$, it follows that we may assume
(\ref{closure}) for every $B_r(y)$ in (\ref{poincare}) when (\ref{eucled2})
holds. Similarly, since the constant $c_0$ in (\ref{poincare}) is
independent of $K$, we may assume as well that every $B_{c_0r}(y)$ in
(\ref{poincare}) has closure in $\Omega$.  

(iii) We can often slightly weaken the assumption in Definition
\ref{poincaredef} that $K$ is an arbitrary compact set in
$\Omega$. For example, in our results where $w(\Omega) < \infty$, it
is generally enough to assume that for each $\epsilon >0$, there is a
particular compact $K$ with $w(\Omega\setminus K)< \epsilon$ such that
(\ref{poincare}) holds. However, in \S 3.4, where we do not assume
$w(\Omega) < \infty$, it is convenient to keep the hypothesis that $K$
is arbitrary.  

\end{remark}

Given a set ${\cal H} \subset Lip_{Q,p}(\Omega)$, define
\bea\label{newset}
\hat{\cal H}= \{ f : \mbox{there exists } \{f^j\} \subset
{\cal H} \mbox{ with  $f^j \to f$  a.e.-$w$} \}.
\eea
It will be useful later to note that if ${\cal H}$ is bounded in
$L^N_w(\Omega)$ for some $N$, then $\hat{{\cal H}}$ is also bounded in
$L^N_w(\Omega)$ by Fatou's lemma; in particular, every $f \in
\hat{\cal H}$ then belongs to $L^N_w(\Omega)$. See (\ref{hat-closure})
for a relationship between $\hat{{\cal H}}$ and the closure of ${\cal
H}$ in $\wa{p}{\Omega}$ in case $w<<\nu$.

We now state our simplest global result. Its proof is given after
Corollary \ref{simplecor*}.

\begin{thm}\label{simpleversion}
Let the assumptions of \S 3.1 hold, $w(\Omega)<\infty$, $1\le p
<\infty$, $1<N \le \infty$ and ${\cal H}\subset
\Lip_{Q,p}(\Omega)$. Suppose that the Poincar\'e property of order $p$
in Definition \ref{poincaredef} holds for ${\cal H}$ and that
\bea\label{sup}
\sup_{f \in {\cal H}} \left\{||f||_{L^N_w(\Omega)} +
||f||_{L^p_\nu(\Omega)} + ||\nabla f||_{{\cal L}^p_\mu(\Omega,Q)}
  \right\}  < \infty.
\eea
Then any sequence $\{f_k\} \subset\hat{{\cal H}}$ has a subsequence
that converges in $L^q_w(\Omega)$ norm for every $1\leq q<N$ to a
function belonging to $L^N_w(\Omega)$.
\end{thm}

Let ${\cal H} \subset Lip_{Q,p}(\Omega)$ and $\hat{\cal H}$ be as in
(\ref{newset}). We reserve the notation $\overline{\cal H}$ for the
closure of ${\cal H}$ in $W^{1,p}_{\nu,\mu}(\Omega, Q)$, i.e., for the
closure of the collection $\{(f, \nabla f): f\in {\cal H}\}$ with
respect to the norm (\ref{pairnormlip}).  Elements of $\overline{\cal
H}$ are viewed as equivalence classes. If $w<<\nu$, then
\bea\label{hat-closure} \{ f : \mbox{there exists $\vec{g}$ such that
} (f,\vec{g})\in \overline{\cal H}\}\subset \hat{\cal H}.
\eea
Indeed, if $(f,\vec{g})\in \overline{\cal H}$, there is a sequence
$\{f^j\} \subset {\cal H}$ such that $(f^j,\nabla
f^j) \ra (f,\vec{g})$ in $\wa{p}{\Omega}$ norm, and consequently $f^j
\ra f$ in $L^p_\nu(\Omega)$. By using a subsequence, we may assume
that $f^j \ra f$ pointwise a.e.-$\nu$, and hence by absolute continuity
that $f^j\to f$ pointwise a.e.-$w$. This proves
(\ref{hat-closure}). In fact, it can be verified by using Egorov's
theorem that
\bea\label{betterhat-closure}
\{f: \text{there exists $\{(f^j,\vec{g^j})\}\subset \overline{\cal H}$
with $f^j \ra f$ a.e.-$w$}\}  \subset \hat{\cal H}.
\eea

Theorem \ref{simpleversion} and (\ref{hat-closure}) immediately imply
the following corollary.

\begin{corollary}\label{simplecor}
Let the assumptions of \S 3.1 hold, $w(\Omega)< \infty$ and
$w<<\nu$. Let $1\le p <\infty$, $1<N \le \infty$, ${\cal H} \subset
Lip_{Q,p}(\Omega)$ and $\overline{\cal H}$ be the closure of ${\cal
H}$ in $W_{\nu,\mu}^{1,p}(\Omega, Q)$. Suppose that the Poincar\'e
property of order $p$ in Definition \ref{poincaredef} holds for ${\cal
H}$ and that
\bea\label{corboundN}
\sup_{f \in {\cal H}} \left\{||f||_{L^N_w(\Omega)} +
||(f,\nabla f)||_{W^{1,p}_{\nu,\mu}(\Omega,Q)}
  \right\}  < \infty.
\eea
Then any sequence $\{f_k\}$ in
\bea
\{ f : \mbox{there exists $\vec{g}$ such that
} (f,\vec{g})\in \overline{\cal H}\} \nonumber
\eea
has a subsequence that converges in $L^q_w(\Omega)$ norm for
$1\le q<N$ to a function that belongs to $L^N_w(\Omega)$.

\end{corollary}

\begin{remark} Corollary \ref{simplecor} may be thought of as
an analogue in the degenerate setting of the
Rellich-Kondrachov theorem since it contains this classical result as
a special case. To see why, set $Q(x) = \mbox{Id}$ and $w=\nu=\mu$
to be Lebesgue measure.  Then, given a bounded sequence
$\{(f_k,\vec{g}_k)\}\subset W^{1,p}_0(\Omega) =
W^{1,p}_{dx,dx,0}(\Omega,Q)$ we may choose $\{f_k^j\}\subset
Lip_0(\Omega)$ with $(f_k^j,\nabla f_k^j)\ra  (f_k,\vec{g}_k)$ in
$W^{1,p}(\Omega)$ norm.  Thus, setting ${\cal H} =
\{f_k^j\}_{k\in\mathbb{N},j>J_k}$ where each $J_k$ is chosen
sufficiently large to preserve boundedness, the classical
Sobolev inequality gives (\ref{corboundN}) with $N=np/(n-p)$ for
$1\leq p<n$.   The Rellich-Kondrachov theorem now follows from
Corollary \ref{simplecor}. 
\end{remark}

We next mention analogues of these results when ${\cal H}$ is replaced
by a set ${\cal W}\subset W^{1,p}_{\nu,\mu}(\Omega,Q)$ with elements
viewed as equivalence classes, assuming that $w <<\nu$. We then
modify Definition \ref{poincaredef} by replacing (\ref{poincare}) with
the analogous estimate
\bea \label{poincare*} \left(\int_{B_r(y)} |f-f_{B_r(y),w}|^{p} dw
\right)^\frac{1}{p} \leq
\epsilon ||(f,\vec{g})||_{W_{\nu,\mu}^{1,p}(B_{c_0r}(y),Q)}\quad
\text{if }(f,\vec{g}) \in {\cal W}.
\eea
The assumption $w<< \nu$ guarantees that the left side of
(\ref{poincare*}) does not change when the first component of a pair
is arbitrarily altered in a set of $\nu$-measure zero.

If Poincar\'e's inequality is known to hold for subsets of
Lipschitz functions in the form (\ref{poincare}), it can often
be extended by approximation to the similar form (\ref{poincare*}) for
subsets of $W^{1,p}_{\nu,\mu}(\Omega,Q)$. Indeed, let us show 
without using weak convergence that if $w<<\nu$ and the Radon-Nikodym
derivative $dw/d\nu  \in L^{p'}_\nu(\Omega), 1/p + 1/p' =1$, then
(\ref{poincare*}) holds with ${\cal W} = W^{1,p}_{\nu,\mu}(\Omega,Q)$
if (\ref{poincare}) holds with ${\cal H} = Lip_{Q,p}(\Omega)$.
This follows easily from Fatou's lemma since if $(f,\vec{g}) \in 
W^{1,p}_{\nu,\mu}(\Omega,Q)$ and we choose $\{f_j\}
\subset Lip_{Q,p}(\Omega)$ with $(f_j, \nabla f_j) \ra (f,\vec{g})$ in
$W^{1,p}_{\nu,\mu}(\Omega,Q)$, then for any ball $B$, since $f_j \ra
f$ in $L^p_\nu(\Omega)$, we have
\[
(f_j)_{B,w} = \frac{1}{w(B)}\int_B f_j\, \frac{dw}{d\nu}\, d\nu \ra
\frac{1}{w(B)}\int_B f\, \frac{dw}{d\nu}\, d\nu = f_{B,w}.
\]
Of course we may also assume that $f_j \ra f$ a.e.-$w$ by selecting a
subsequence of $\{f_j\}$ which converges to $f$ a.e.-$\nu$.
The same argument shows that if (\ref{poincare*}) holds for all pairs
in any set ${\cal W}\subset W^{1,p}_{\nu,\mu}(\Omega, Q)$, then it
also holds for pairs in the closure $\overline{\cal W}$ of ${\cal W}$
in $W^{1,p}_{\nu,\mu}(\Omega,Q)$. Moreover, if all balls $B$ in
question satisfy $\overline{B} \subset \Omega$ (cf. (\ref{closure})),
then the assumption can clearly be weakened to $dw/d\nu \in
L^{p'}_{\nu, loc}(\Omega)$. As we observed in Remark
\ref{various}(ii), the balls in (\ref{poincare}) can be assumed to
satisfy (\ref{closure}) provided (\ref{eucled2}) is true.

Analogues of Theorem \ref{simpleversion} and Corollary
\ref{simplecor} for a set ${\cal W}\subset
W^{1,p}_{\nu,\mu}(\Omega,Q)$ are given in the next result, which
also includes the Rellich-Kondrachov theorem as a special case.

\begin{thm}\label{simpleversion*}
Let the assumptions of \S 3.1 hold, $w(\Omega)<\infty$ and
$w<<\nu$. Let $1\le p <\infty$, $1<N \le \infty$ and ${\cal W}
\subset W^{1,p}_{\nu,\mu}(\Omega,Q)$. Suppose that the Poincar\'e
property in Definition \ref{poincaredef} holds, but in the modified
form given in (\ref{poincare*}), and that
\bea\label{sup*}
\sup_{(f,\vec{g}) \in {\cal W}} \left\{||f||_{L^N_w(\Omega)} +
||(f,\vec{g})||_{W^{1,p}_{\nu,\mu}(\Omega,Q)}  \right\}  < \infty.
\eea
Let
\[
\hat{\cal W} = \{f: \text{there exists $\{(f^j,\vec{g^j})\} \subset
{\cal W}$ with $f^j \ra f$ a.e.$-w$}\}.
\]
Then any sequence in $\hat{{\cal W}}$ has a subsequence that converges
in $L^q_w(\Omega)$ norm for every $1\leq q<N$ to a function belonging
to $L^N_w(\Omega)$. In particular, if $\overline{\cal W}$ denotes the
closure of ${\cal W}$ in $W_{\nu,\mu}^{1,p}(\Omega, Q)$, then the same
is true for any sequence in
\bea
\{ f : \mbox{there exists $\vec{g}$ such that
} (f,\vec{g})\in \overline{\cal W}\}. \nonumber
\eea
\end{thm}

As a corollary, we obtain a result for arbitrary sequences $\{(f_k,
\vec{g_k})\}$ which are bounded in $W^{1,p}_{\nu,\mu}(\Omega,Q)$ and
whose first components $\{f_k\}$ are bounded in $L^N_w(\Omega)$.

\begin{corollary}\label{simplecor*}
Let the assumptions of \S 3.1 hold, $w(\Omega)<\infty$, $w<<\nu$,
$1\le p <\infty$ and $1<N \le \infty$. Suppose that the Poincar\'e
property in Definition \ref{poincaredef} holds for all of
$W^{1,p}_{\nu,\mu}(\Omega,Q)$, i.e., Definition \ref{poincaredef}
holds with (\ref{poincare}) replaced by (\ref{poincare*}) for ${\cal W} =
W^{1,p}_{\nu,\mu}(\Omega,Q)$. Then if $\{(f_k, \vec{g_k})\}$ is any
sequence in $W_{\nu,\mu}^{1,p}(\Omega, Q)$ such that
\bea
\sup_k \left[ ||f_k||_{L^N_w(\Omega)} + ||(f_k,\vec{g_k})||_{W^{1,
p}_{\nu,\mu}(\Omega,Q)}\right] < \infty,
\nonumber
\eea
there is a subsequence of $\{f_k\}$ that converges in $L^q_w(\Omega)$
norm for $1\le q<N$ to a function belonging to $L^N_w(\Omega)$. If in
addition $dw/d\nu \in L^{p'}_\nu(\Omega), 1/p + 1/p' =1$, the
conclusion remains valid if the Poincar\'e property holds just for
$Lip_{Q,p}(\Omega)$.
\end{corollary}

In fact, the first conclusion in Corollary \ref{simplecor*} follows
by applying Theorem \ref{simpleversion*} with ${\cal W}$ chosen to be
the specific sequence $\{(f_k, \vec{g_k})\}_k$ in question, and the
second statement follows from the first one and our observation above
that (\ref{poincare*}) holds with ${\cal W} =
W^{1,p}_{\nu,\mu}(\Omega,Q)$ if $dw/d\nu \in L^{p'}_\nu(\Omega), 1/p + 1/p'
=1$, and if (\ref{poincare}) holds with ${\cal H} = Lip_{Q,p}(\Omega)$. \\

\noindent {\bf Proofs of Theorems \ref{simpleversion} and
\ref{simpleversion*}.} We will concentrate on the proof of Theorem
\ref{simpleversion}. The proof of Theorem \ref{simpleversion*} is
similar and omitted. We begin with a useful covering lemma.

\begin{lemma}\label{help} Let the assumptions of \S 3.1 hold and
$w(\Omega)< \infty$. Fix $p\in [1,\infty)$ and a set ${\cal H}
\subset Lip_{Q,p}(\Omega)$. Suppose the Poincar\'e property of order
$p$ in Definition \ref{poincaredef} holds for ${\cal H}$, and let
$\kappa$ be as in (\ref{tri}) and $c_0$ be as in
(\ref{poincare}). Then for every $\epsilon>0$, there are positive
constants $r=r(\epsilon,\kappa, c_0), M= M(\kappa, c_0)$ and a
finite collection $\{B_{r}(y_k)\}_k$ of $d$-balls, so that
\bea \label{lemma_i}&& (i) \quad
w\big(\Omega\setminus \disp\bigcup_{k}
B_{r}(y_k)\big)<\epsilon ,\\
\label{lemma_ii}&&  (ii) \quad\sum_{k} \chi_{B_{c_0r}(y_k)}(x)\leq M
\quad \text{for all $x\in \Omega$}, \\
\label{lemma_iii}&& (iii) \quad ||f-f_{B_{r}(y_k),w}||_{L^{p}_w(B_{r}(y_k))}
\leq \epsilon ||(f,\nabla f)||_{W^{1,p}_{\nu,\mu}(B_{c_0r}(y_k),Q)}
\eea
for all $f \in {\cal H}$ and all $k$. Note that $M$ is
independent of $\epsilon$.

\end{lemma}

\noindent {\bf Proof:}  We first recall the ``swallowing'' property of
$d$-balls: There is a constant $\gamma \ge 1$ depending only
on $\kappa$ so that if $x, y \in \Omega$, $0< r_1\leq r_2< \infty$ and
$B_{r_1}(x)\cap B_{r_2}(y)\neq\emptyset$, then
\bea \label{swallowing}B_{r_1}(x)\subset B_{{\gamma}r_2}(y).
\eea
Indeed, by \cite[Observation 2.1]{CW1}, $\gamma$ can be chosen to be
$\kappa+2\kappa^2$. 

Fix $\epsilon>0$. Since $w(\Omega)< \infty$, there is a compact set
$K\subset \Omega$ with $w(\Omega\setminus K) <\epsilon$. Let $\delta'
= \delta'(\epsilon)$ be as in Definition \ref{doublingdef} for $K$,
and let $\delta = \delta(\epsilon)$ be as in (\ref{poincare}). Fix $r$
with $0<r<\textrm{min}\{\delta, \delta'/(c_0\gamma)\}$ where $c_0$ is
as in (\ref{poincare}). For each $x\in K$, use (\ref{new2}) to pick
$s(x,r)>0$ so that $D_{s(x,r)}(x)\subset B_{r/\gamma}(x)$. Since $K$
is compact, there are finitely many points $\{x_j\}$ in $K$ so that
$K\subset \cup_j B_{r/\gamma}(x_j)$.  Choose a maximal pairwise
disjoint subcollection $\{B_{r/\gamma}(y_k)\}$ of
$\{B_{r/\gamma}({x_j})\}$. We will show that the collection
$\{B_r(y_k)\}$ satisfies (\ref{lemma_i})--(\ref{lemma_iii}).

To verify (\ref{lemma_i}), it is enough to show that $K\subset \cup_k
B_r(y_k)$. Let $y\in K$.  Then $y\in B_{r/\gamma}(x_j)$ for some
$x_j$. If $x_j = y_k$ for some $y_k$ then $y\in B_r(y_k)$.  If
$x_j\neq y_k$ for all $y_k,$ there exists $y_\ell$ so that
$B_{r/\gamma}(y_\ell) \cap B_{r/\gamma}(x_j)\neq\emptyset$.  Then
$B_{r/\gamma}(x_j)\subset B_{r}(y_\ell)$ by (\ref{swallowing}), and so
$y\in B_r(y_\ell)$.  In either case, we obtain $y\in \disp\cup_k
B_r(y_k)$ as desired. 

To verify (\ref{lemma_ii}), suppose that $\{k_i\}_{i=1}^L$ satisfies
$\cap_{i=1}^L B_{c_0r}(y_{k_i}) \neq \emptyset$.
 Then by
(\ref{swallowing}), $B_{c_0r}(y_{k_i}) \subset
B_{c_0{\gamma}r}(y_{k_1})$ for $1\leq i\leq L$. Since $\gamma, c_0 \ge
1$, we have $B_{r/\gamma}(y_k) \subset B_{c_0r}(y_k)$ for all $k$, and
consequently
$$
\cup B_{r/\gamma}(y_{k_i}) \subset \cup B_{c_0r}(y_{k_i}) \subset
B_{c_0\gamma r}(y_{k_1}).
$$
By construction, $\{B_{r/\gamma}(y_k)\}$ is pairwise disjoint in
$k$. Since $0<r/\gamma< c_0\gamma r<\delta'$, the corresponding
constant ${\cal C}$ in the definition of geometric doubling depends
only on $(c_0\gamma r)/(r/\gamma) = c_0\gamma^2$, i.e., ${\cal C}$
depends only on $\kappa$ and $c_0$. Choosing $M$ to be this constant,
we obtain that $L\leq M$ as desired. The same argument shows that the
collection $\{B_{c_0r}(y_k)\}$ has the stronger bounded intercept
property with the same bound $M$, i.e., any ball in the collection
intersects at most $M-1$ others.

Finally, let us verify (\ref{lemma_iii}).  Recall that
$0<r<\delta$ by construction. Hence (\ref{poincare}) implies that for
each $k$ and all $f \in {\cal H }$,
\begin{equation} ||f-f_{B_{r}(y_k),w}||_{L^{p}_w(B_{r}(y_k))}
\leq \epsilon ||(f,\nabla f)||_{W^{1,p}_{\nu, \mu}(B_{c_0r}(y_k),Q)},
\end{equation}
as required. This completes the proof of Lemma \ref{help}. $\Box$

The proof of Theorem \ref{simpleversion} will be deduced from Theorem
\ref{general} by choosing $\mathscr{X}(\Omega) =
L^p_\nu(\Omega) \times {\cal L}^{p}_\mu(\Omega,Q)$ and considering the
product space
\[
{\cal B}_{N,{\mathscr{X}}(\Omega)} =
L^N_w(\Omega)\times \left(L^p_\nu(\Omega) \times {\cal
L}^{p}_\mu(\Omega,Q)\right).
\]
We always choose $\Sigma$ to be the Lebesgue measurable subsets of
$\Omega$ and $\Sigma_0 = \{B_r(x): r>0,x\in\Omega\}$. Note that
$\mathscr{X}(\Omega)$ and ${\cal B}_{N,{\mathscr{X}}(\Omega)}$ are
normed linear spaces (even Banach spaces), and the norm in ${\cal
B}_{N,{\mathscr{X}}(\Omega)}$ is
\bea ||(h, (f,\vec{g}))||_{{\cal B}_{N,{\mathscr{X}}(\Omega)}} =
||h||_{L^N_w(\Omega)} + ||f||_{L^p_\nu(\Omega)} +
||\vec{g}||_{{\cal L}^{p}_\mu(\Omega,Q)}.
\eea
The roles played in \S 1 by $\bf{g}$ and $(f,\bf{g})$ are now played
by $(f,\vec{g})$ and $(h, (f,\vec{g}))$ respectively.

Let us verify properties (A) and ($B_p$) in \S 1 with
${\mathscr{X}}(\Omega)$ and $\Sigma_0$ chosen as above.  To verify
(A), fix $B\in \Sigma_0$ and $(f,\vec{g})\in
{\mathscr{X}}(\Omega)$. Clearly $f\chi_B \in L^p_\nu(\Omega)$ since $f
\in L^p_\nu(\Omega)$. Also,
\bea \int_\Omega \Big((\vec{g}\chi_{B})'Q (\vec{g} \chi_{B})
\Big)^\frac{p}{2}d\mu &=&  \int_{B} \Big(\vec{g}\,'Q(x)\vec{g}
\Big)^\frac{p}{2}d\mu \nonumber\\
&\leq& \int_\Omega \Big(\vec{g}\,'Q(x)\vec{g} \Big)^\frac{p}{2}d\mu
<\infty.\nonumber
\eea
Thus $(f,\vec{g})\chi_{B}\in {\mathscr{X}}(\Omega)$ and property (A)
is proved.

To verify ($B_p$), let $\{B_l\}$ be a finite collection
of $d$-balls satisfying $\sum_{l} \chi_{B_l}(x)\leq C_1$ for all $x
\in\Omega$. Then if $(f,\vec{g})\in {\mathscr{X}}(\Omega)$,
$$
\disp\sum_{l}||(f,\vec{g}) \chi_{B_l}||^p_{\mathscr{X}(\Omega)}
=  \disp\sum_l \left(||f\chi_{B_l}||_{L^p_\nu(\Omega)} +
||\vec{g} \chi_{B_l}||_{{\cal L}^{p}_\mu(\Omega,Q)}\right)^p
$$
$$
\le 2^{p-1} \disp\sum_l \left(||f\chi_{B_l}||^p_{L^p_\nu(\Omega)} +
||\vec{g} \chi_{B_l}||^p_{{\cal L}^{p}_\mu(\Omega,Q)}\right)
$$
$$
= 2^{p-1} \int_\Omega |f|^p\left(\disp\sum_{l}\chi_{B_l}\right) d\nu +
\int_\Omega \left(\vec{g}\,'Q\vec{g} \right)^{\frac{p}{2}} \left(\disp
\sum_{l}\chi_{B_l}\right) d\mu
$$
$$
\le 2^{p-1}C_1 \left(||f||^p_{L^p_\nu(\Omega)} +
||\vec{g}||^p_{{\cal L}^{p}_\mu(\Omega,Q)}\right) \le 2^pC_1
  ||(f, \vec{g})||^p_{{\mathscr{X}}(\Omega)}.
$$
This verifies ($B_p$) with $C_2$ chosen to be $2^pC_1$.

The proof of Theorem \ref{simpleversion} is now very simple. Let
${\cal H}$ satisfy its hypotheses and choose ${\cal S}$ in Theorem
\ref{general} to be the set
\[
{\cal S} =  \left\{(f, (f,\nabla f)) : f \in {\cal H}\right\}.
\]
Note that ${\cal S}$ is a bounded subset of ${\cal
B}_{N,{\mathscr{X}}(\Omega)}$ by hypothesis (\ref{sup}). Next, in
order to choose the pairs $\{E_\ell, F_\ell\}_\ell$ and verify
conditions (i)--(iii) of Theorem \ref{general} (see (\ref{bddov}) and
(\ref{pc})), we appeal to Lemma \ref{help}. Given $\epsilon>0$, let
$\{E_\ell, F_\ell\}_\ell =  \{B_r(y_k),B_{c_0r}(y_k)\}_{k}$ where
$\{y_k\}$ and $r$ are as in Lemma \ref{help}. Then $E_\ell, F_\ell\in
\Sigma_0$, and conditions (i)--(iii) of Theorem \ref{general} are
guaranteed by Lemma \ref{help}. Finally, by noting that the set
$\hat{\cal H}$ defined in (\ref{newset}) is the same as the set
$\hat{\cal S}$ defined in (\ref{hatS}), the conclusion of Theorem
\ref{simpleversion} follows from Theorem \ref{general}.  $\Box$\\

For special domains $\Omega$ and special choices of $N$, the
boundedness assumption (\ref{sup}) (or (\ref{corboundN})) can
be weakened to
\bea\label{weaksup}
\sup_{f \in {\cal H}} \left\{||f||_{L^p_\nu(\Omega)} +
||\nabla f||_{{\cal L}^p_\mu(\Omega,Q)}  \right\} =
\sup_{f \in {\cal H}} ||(f,\nabla
f)||_{W^{1,p}_{\nu,\mu}(\Omega,Q)}   < \infty.
\eea
This is clearly the case for any $\Omega$ and $N$ for which there
exists a global Sobolev-Poincar\'e estimate that bounds
$||f||_{L^N_w(\Omega)}$ by $||(f,\nabla
f)||_{W^{1,p}_{\nu,\mu}(\Omega, Q)}$ for all $f \in {\cal
H}$.  We now formalize this situation assuming that $w <<\nu$. In the 
appendix, we consider a case when $w<<\nu$ fails. 

The form of the global Sobolev-Poincar\'e estimate we will use is
given in the next definition.  It  guarantees that (\ref{sup}) and
(\ref{weaksup}) are the same when $N=p\sigma$.

\begin{defn}\label{strongsobdef}  Let $1\leq p<\infty$ and ${\cal H}
\subset Lip_{Q,p}(\Omega)$.  Then the \emph{global Sobolev property of
order $p$ holds for} ${\cal H}$ if there are constants $C>0$ and
$\sigma>1$ so that
\bea\label{strongsob}
||f||_{L^{p\sigma}_w(\Omega)} \leq C ||(f,\nabla f)||_{\wa{p}{\Omega}}
  \quad\text{for all $f \in {\cal H}$}.
\eea
\end{defn}

If $w<<\nu$, then (\ref{strongsob}) extends to $(f,\vec{g}) \in
\overline{\cal H}$. In fact, let $(f,\vec{g}) \in \overline{{\cal H}}$
and choose $\{f_j\}\subset {\cal H}$ with $(f_j, \nabla f_j)\to
(f,\vec{g})$ in $W^{1,p}_{\nu,\mu}(\Omega,Q)$. Then $f_j \ra f$ in
$L^p_{\nu}(\Omega)$ norm, and by choosing a subsequence we may assume
that $f_j \ra f$ a.e.-$\nu$. Hence $f_j \ra f$ a.e.-$w$ because
$w<<\nu$. Since each $f_j$ satisfies (\ref{strongsob}), it follows that
\bea\label{strongsob**}
||f||_{L^{p\sigma}_w(\Omega)} \le C ||(f,\vec{g})||_{\wa{p}{\Omega}}
\quad\text{if $(f,\vec{g})\in \overline{\cal H}$}.
\eea
Under the same assumptions, namely that Definition \ref{strongsobdef}
holds for a set ${\cal H} \subset Lip_{Q,p}(\Omega)$ and that
$w<<\nu$, the same sequence $\{f_j\}$ as above is also bounded in
$L_w^{p\sigma}(\Omega)$ norm and so satisfies $(f_j)_{E,w} \ra f_{E,w}$
for measurable $E$ by the same weak convergence argument given
after the statement of Theorem \ref{general}. Hence the
Poincar\'e estimate in Definition \ref{poincaredef} also extends to
$\overline{\cal H}$ in the same form as (\ref{poincare*}), with 
${\cal W}$ there replaced by $\overline{\cal H}$, i.e., 
\bea \label{poincare**} \left(\int_{B_r(y)} |f-f_{B_r(y),w}|^{p} dw
\right)^\frac{1}{p} \leq \epsilon
||(f,\vec{g})||_{W_{\nu,\mu}^{1,p}(B_{c_0r}(y),Q)}\quad 
\text{if }(f,\vec{g}) \in \overline{\cal H}.
\eea
Hence, we immediately obtain the next result by choosing ${\cal W}
=\overline{\cal H}$ and $N=p\sigma$ in Theorem \ref{simpleversion*}.

\begin{thm} \label{appth1}   Let the assumptions of \S 3.1 hold, 
$w(\Omega)< \infty$ and $w<<\nu$. Fix $p\in [1,\infty)$ and a set
${\cal H} \subset Lip_{Q,p}(\Omega)$.  Suppose the Poincar\'e and 
global Sobolev properties of order $p$ in Definitions
\ref{poincaredef} and \ref{strongsobdef} hold for ${\cal H}$, and let
$\sigma$ be as in (\ref{strongsob}). If $\{(f_k,\vec{g_k})\}$ is a
sequence in $\overline{{\cal H}}$ with
\bea \label{w1pbound}\sup_k
||(f_k,\vec{g_k})||_{W^{1,p}_{\nu,\mu}(\Omega,Q)} < \infty,
\eea
then $\{ f_k\}$ has a subsequence which converges in
$L^q_w(\Omega)$ for $1 \le q< p\sigma$, and the limit of the
subsequence belongs to $L^{p\sigma}_w(\Omega)$.
\end{thm}

A result for the entire space $W^{1,p}_{\nu,\mu}(\Omega, Q)$ follows
by choosing ${\cal H}= Lip_{Q,p}(\Omega)$ in Theorem \ref{appth1} or
Corollary \ref{simplecor}: 

\begin{corollary}\label{appcor1} Suppose that the hypotheses of
Theorem \ref{appth1} hold for ${\cal H} = Lip_{Q,p}(\Omega)$.
If $\{(f_k,\vec{g_k})\} \subset \wa{p}{\Omega}$ and (\ref{w1pbound})
is true then $\{ f_k\}$ has a subsequence which converges in
$L^q_w(\Omega)$ for $1 \le q< p\sigma$, and the limit of the
subsequence belongs to $L^{p\sigma}_w(\Omega)$.
\end{corollary}

See the Appendix for analogues of Theorem \ref{appth1} and Corollary
\ref{appcor1} without the assumption $w<<\nu$.

\subsection{Local Compactness Results for Degenerate Spaces}

In this section, for general bounded measurable sets $\Omega'$ with
$\overline{\Omega'}\subset \Omega$, we study compact embedding of
subsets of $W^{1,p}_{\nu,\mu}(\Omega,Q)$ into $L^q_w(\Omega')$
without assuming a global Sobolev estimate for $\Omega$ or $\Omega'$
and without assuming $w(\Omega) <\infty$.  For some applications, see
the comment at the end of the section.

The theorems below will assume a much weaker condition than the global
Sobolev estimate (\ref{strongsob}), namely the following local
estimate.

\begin{defn}\label{localsobdef} Let $1\le p<\infty$. We say that the
\emph{local Sobolev property of order $p$} holds if for some fixed
constant $\sigma >1$ and every compact set $K\subset\Omega$, there is
a constant $r_1>0$ so that for all $d$-balls $B= B_r(y)$ with $y\in K$
and $0<r<r_1$,
\bea \label{sob}
||f||_{L^{p\sigma}_w(B)} \le C(B)\,
||(f,\nabla f)||_{W^{1,p}_{\nu,\mu}(\Omega, Q)} \quad \text{if } f\in
Lip_0(B),
\eea
where $C(B)$ is a positive constant independent of $f$. We will
view any $f\in Lip_0(B)$ as extended by $0$ to all of $\Omega$.
\end{defn}

\begin{remark}\label{lots}
(i) A more standard assumption than (\ref{sob}) is a normalized
inequality that includes a factor $r$ in the gradient term on the right
side:
\bea\label{sob2} \left(\frac{1}{w(B_r(y))}\int_{B_r(y)}
|f|^{{p}\sigma}dw\right)^\frac{1}{{p}\sigma} \leq
C\left(\frac{1}{\nu(B_r(y))} \int_{B_r(y)} |f|^{p}d\nu\right)^\frac{1}{p}
\nonumber\\ + Cr \left(\frac{1}{\mu(B_r(y))}\int_{B_r(y)} |\sqrt{Q}\nabla
f|^{p}d\mu\right)^\frac{1}{{p}},
\eea
with $C$ independent of $r, y$;  see e.g. \cite{SW1} and \cite{R1} in
the unweighted case with $p=2$. Clearly (\ref{sob2}) is a stronger
requirement than (\ref{sob}).

(ii) In the classical $n$-dimensional elliptic case for linear second
order equations in divergence form, $Q$ satisfies $c|\xi|^2 \le
Q(x,\xi) \le C|\xi|^2$ for some fixed constants $c, C>0$ and $d$ is the
standard Euclidean metric $d(x,y) = |x-y|$. For $1\le p<n$ and
$\sigma = n/(n-p)$, (\ref{sob}) then holds with $dw= d\nu =d\mu = dx$
since the corresponding version of (\ref{sob2}) is true with
$|\sqrt{Q}\nabla f|$ replaced by $|\nabla f|$.
\end{remark}

We will also use a notion of Lipschitz cutoff functions on $d$-balls:

\begin{defn}\label{cutoff}  For $s\geq 1$, we say that the
\emph{cutoff property of order $s$} holds for $\mu$ if for each
compact $K\subset\Omega$, there exists $\delta =\delta(K)>0$ so that
for every d-ball $B_r(y)$ with $y\in K$ and $0<r<\delta,$ there is a
function $\phi\in Lip_0(\Omega)$ and a constant $\gamma=\gamma(y,r) \in
(0,r)$ satisfying

\medskip

\noindent (i) $\;$ $0\leq \phi \leq 1$ in $\Omega$,  \\
\noindent (ii) $\;$ supp $\phi \subset B_r(y)$ and $\phi =1$ in
$B_\gamma(y)$,  \\
\noindent (iii) $\;$ $\nabla \phi \in {\cal L}^s_\mu(\Omega,Q)$.
\end{defn}

Since $\mu$ is always assumed to be locally finite, the strongest form
of Definition \ref{cutoff}, namely the version with $s=\infty$,
automatically holds if $Q$ is locally bounded in $\Omega$ and
(\ref{eucled2}) is true; recall that we always assume (\ref{new2}). To
see why, fix a compact set $K\subset \Omega$ and consider $B_r(y)$
with $y\in K$ and $r<1$.  Use (\ref{new2}) to choose open Euclidean
balls $D', D$ with common center $y$ such that $\overline{D'} \subset
D \subset B_r(y) (\subset \Omega \text{ by definition})$. Construct a
smooth function $\phi$ in $\Omega$ with support in $D$ such that $0\le
\phi \le 1$ and $\phi = 1$ on $D'$. By (\ref{eucled2}), there is
$\gamma >0$ such that $B_\gamma(y)\subset D'$. Then $\phi$ satisfies
parts (i)-(iii) of Definition \ref{cutoff} with $s=\infty$; for (iii),
we use the fact that $\nabla\phi$ has compact support in $\Omega$
together with local boundedness of $Q$ and local finiteness of $\mu$.

To compensate for the lack of a global Sobolev estimate, given ${\cal
H} \subset Lip_{Q,p}(\Omega)$, we will assume in conjunction with the
cutoff property of some order $s\ge p\sigma'$ that for every compact
set $K\subset \Omega$, there exists $\delta = \delta(K) >0$ such that
for every $d$-ball $B$ with center in $K$ and radius less than
$\delta$, there is a constant $C_1(B)$ so that
\bea\label{globpc}
||f||_{L^{pt'}_\mu(B)} \leq C_1(B)\, ||(f, \nabla
  f)||_{W^{1,p}_{\nu,\mu}(\Omega,Q)} \quad \text{if } f\in {\cal H},
\eea
where  $t=s/p$ and $1/t+ 1/t'=1$. Note that $1\leq t'\leq \sigma$
since $s\ge p\sigma'$.

\begin{remark} \label{extracondition} Inequality (\ref{globpc}) is
different in nature from (\ref{sob}) even if $t'=\sigma$ and $w=\mu$
since there is a restriction on supports in (\ref{sob}) but not in
(\ref{globpc}). However, (\ref{globpc}) implies (\ref{sob}) when
$s=p\sigma'$, $w=\mu$ and ${\cal H}$ contains all Lipschitz functions
with support in any ball. On the other hand, (\ref{globpc}) is often
automatic if $\mu=\nu$. For example, as mentioned earlier, if  $Q$ is
locally bounded and (\ref{eucled2}) is true, then the cutoff property
holds with $s=\infty$, giving $t= \infty$ and $t'=1$. In this case,
when $\mu = \nu$, the left side of (\ref{globpc}) is clearly smaller
than the right side (in fact smaller than $||f||_{L^p_\nu(\Omega)}$).

\end{remark}

We can now state our main local result.

\begin{thm}\label{newappth2} Let the assumptions of \S 3.1 and
condition (\ref{eucled2}) hold, and let $w<<\nu$. Fix $p\in
[1,\infty)$ and suppose the Poincar\'e property of order $p$ in
Definition \ref{poincaredef} holds for a fixed set ${\cal H}\subset
Lip_{Q,p}(\Omega)$ and the local Sobolev property of order $p$ in
Definition \ref{localsobdef} holds. Assume the cutoff
property of some order $s\geq p\sigma'$ is true for $\mu$, with
$\sigma$ as in (\ref{sob}), and that (\ref{globpc}) holds for ${\cal
H}$ with $t=s/p$. Then for every $\{(f_k,\vec{g_k})\} \subset
\overline{\cal H}$ that is bounded in $W^{1,p}_{\nu,\mu}(\Omega,Q)$
norm, there is a subsequence $\{f_{k_i}\}$ of $\{f_k\}$ and an
$f \in L^{p\sigma}_{w,loc}(\Omega)$ such that $f_{k_i} \ra f$ pointwise
a.e.-$w$ in $\Omega$ and in $L^q_w(\Omega')$ norm for all $1\le q
<p\sigma$ and every bounded measurable $\Omega'$ with
$\overline{\Omega'} \subset \Omega$.
\end{thm}

See the Appendix for a version of Theorem \ref{newappth2} without
assuming $w<<\nu$.

Recall that $\overline{\cal H} = W^{1,p}_{\nu,\mu}(\Omega,Q)$ if
${\cal H} = Lip_{Q,p}(\Omega)$. In the important case when $Q \in
L^\infty_{loc}(\Omega)$, Theorem \ref{newappth2} and Remark
\ref{extracondition} immediately imply the next result.

\begin{corollary}\label{newappcor2} Let $Q$ be locally bounded in
$\Omega$ and suppose that (\ref{eucled2}) holds.  Fix $p\in
[1,\infty)$, and with $w=\nu=\mu$, assume the Poincar\'e property of
order $p$ holds for $Lip_{Q,p}(\Omega)$ and the local Sobolev property
of order $p$ holds. Then for every bounded sequence $\{(f_k,\vec{g_k})\}
\subset W^{1,p}_{w,w}(\Omega,Q)$, there is a subsequence $\{f_{k_i}\}$
of $\{f_k\}$ and a function $f \in L^{p\sigma}_{w,loc}(\Omega)$ such
that $f_{k_i} \ra f$ pointwise a.e.-$w$ in $\Omega$ and in
$L^q_w(\Omega')$ norm, $1\le q <p\sigma$, for every bounded measurable
$\Omega'$ with $\overline{\Omega'} \subset \Omega$.

\end{corollary}

\medskip

\noindent{\bf Proof of Theorem \ref{newappth2}:}  We begin by using the
cutoff property in Definition \ref{cutoff} to construct a partition of
unity relative to $d$-balls and compact subsets of $\Omega$.

\begin{lemma}\label{partitionofunity} Fix $\Omega$ and $s\geq 1$, and
suppose the cutoff property of order $s$ holds for $\mu$.  If $K$ is a
compact subset of $\Omega$ and $r>0$, there is a finite collection of
$d$-balls $\{B_r(y_j)\}$ with $y_j \in K$ together with Lipschitz
functions $\{\psi_j\}$ on $\Omega$ such that $supp \,\psi_j \subset
B_r(y_j)$ and

\noindent (a) $K\subset \disp\bigcup_j B_r(y_j)$,\\
\noindent (b) $0\leq \psi_j\leq 1$ in $\Omega$ for each $j$, and
$\disp\sum_j \psi_j(x) =1$ for all $x\in K$,\\
\noindent (c) $\nabla \psi_j \in {\cal L}^s_\mu(\Omega,Q)$ for each $j$.
\end{lemma}

\noindent {\bf Proof:} The argument is an adaptation of one in
\cite{Ru} for the usual Euclidean case. The authors thank
D. D. Monticelli for related discussions.  Fix $r>0$ and a compact set
$K\subset\Omega$, and set $\beta = \min\{\delta/2,r\}$ for $\delta=
\delta(K)$ as in Definition \ref{cutoff}. Since $\beta< \delta$,
Definition \ref{cutoff} implies that for each $y\in K$, there exist
$\gamma(y) \in (0,\beta)$ and $\phi_y(x)\in Lip(\Omega)$ so that
$0\le \phi_y\le 1$ in $\Omega$, $supp \,\phi_y \subset B_{\beta}(y))$,
$\phi_y = 1$ in $B_{\gamma(y)}(y)$ and $\nabla\phi_y \in {\cal
L}^s_\mu(\Omega,Q)$. The collection $\{B_{\gamma(y)}(y)\}_{y\in K}$
covers $K$, so by (\ref{new2}) and the compactness of $K$, there is a
finite subcollection $\{B_{\gamma(y_j)}(y_j)\}_{j=1}^m$ whose union
covers $K$.  Part (a) follows since $\gamma(y_j)<r$.
Next let $\phi_j(x)=\phi_{y_j}(x)$ and define $\{\psi_j\}_{j=1}^m$ as
follows: set $\psi_1 = \phi_1$ and $\psi_j=(1-\phi_1)\cdots
(1-\phi_{j-1})\phi_j$ for $j=2,..,m$.  Then each $\psi_j$ is a
Lipschitz function in $\Omega$, and $supp \,\phi_j \subset B_r(y_j)$
since $\beta <r$. Also, $0\leq \psi_j\leq 1$ in $\Omega$ and
\bea \disp\sum_{j=1}^m \psi_j(x)= 1-\prod_{j=1}^m (1-\phi_j(x)), \quad
x\in \Omega.\nonumber
\eea
If $x\in K$ then $x\in B_{\gamma(y_j)}(y_j)$ for some $j$. Hence some
$\phi_j(x)=1$ and consequently $\sum_j \psi_j(x)= 1$ . This
proves part (b). Lastly, we use Leibniz's product rule to compute
$\nabla \psi_j$  and then apply Minkowski's inequality $j$ times to
obtain part (c) from the fact that $\nabla\phi_j\in {\cal
L}^s_\mu(\Omega,Q)$.  $\Box$\\

The next lemma shows how the local Sobolev estimate (\ref{sob}) and
Lemma \ref{partitionofunity} lead to a local analogue of the global
Sobolev estimate (\ref{strongsob}).

\begin{lemma}\label{localglobal} Let $\Omega'$ be a bounded measurable
set with $\overline{\Omega'} \subset \Omega$. Suppose that both Definition
\ref{localsobdef} and the cutoff property for $\mu$ of some order
$s\ge p\sigma'$ hold, and also that (\ref{globpc}) holds with $t=s/p$
for a fixed set ${\cal H} \subset Lip_{loc}(\Omega)$. Then there is a finite
constant $C(\Omega')$ such that
\bea\label{localglobalsob}
||f||_{L^{p\sigma}_w(\Omega')} \le C(\Omega')\, ||(f,\nabla
  f)||_{W^{1,p}_{\nu,\mu}(\Omega,Q)}\quad \text{if } f \in {\cal H}.
\eea
\end{lemma}

\noindent{\bf Proof:}  Let $r_1$ be as in Definition \ref{localsobdef}
relative to the compact set $\overline{\Omega'}\subset \Omega$, and
let $\delta$ be as in (\ref{globpc}). Use Lemma \ref{partitionofunity}
to cover $\overline{\Omega'}$ by the union of a finite number of $d$-balls
$\{B_j\}$ each of radius smaller than $\min\{r_1,\delta\}$. Associated
with this cover is a collection $\{\psi_j\}\subset Lip(\Omega)$ with
$supp \,\psi_j\subset B_j$, $\sum_j\psi_j=1$ in $\Omega'$, and
$\nabla\psi_j\in {\cal L}^s_\mu(\Omega,Q)$. If $f \in {\cal H}$, then
\begin{eqnarray}\label{3.6-1} ||f||_{L^{p\sigma}_w(\Omega')} =
||f\sum_j \psi_j ||_{L^{p\sigma}_w(\Omega')} \leq \disp\sum_j ||\psi_j
f||_{L^{p\sigma}_w(B_j)}.
\eea
Since $\psi_j f\in Lip_0(B_j)$, (\ref{sob}) and the product rule give
\begin{eqnarray}\label{3.6-2}
&&  ||\psi_j f||_{L^{p\sigma}_w(B_j)}\le C(B_j)\,||(\psi_j f,
\nabla(\psi_jf))||_{W^{1,p}_{\nu,\mu}(B_j,Q)}  \nonumber\\
&=& C(B_j)\left(||\psi_j f||_{L^{p}_\nu(B_j)} + ||\sqrt{Q}\nabla(\psi_j
f)||_{L^p_\mu(B_j)}\right)   \nonumber\\
&\leq& C(B_j)\left(||\psi_j f||_{L^p_\nu(B_j)} + ||\psi_j\sqrt{Q}\nabla
f||_{L^p_\mu(B_j)} + ||f\sqrt{Q}\nabla \psi_j||_{L^p_\mu(B_j)} \right)
\nonumber\\ &\leq& C(B_j)\left(||(f,\nabla f)||_{\wa{p}{\Omega}} +
||f\sqrt{Q}\nabla \psi_j||_{L^p_\mu(B_j)}\right),
\eea
where we have used $|\psi_j| \le 1$. We will estimate the second term
on the right of (\ref{3.6-2}) by using (\ref{globpc}). Recall that
$t=s/p \ge \sigma'$ and  $1/t + 1/t'=1$. Let
\[
\overline{C} = \max_j ||\sqrt{Q}\nabla\psi_j||_{L_\mu^s(B_j)}.
\]
By H\"older's inequality and (\ref{globpc}),
\bea \label{3.6-3}
||f\sqrt{Q}\nabla \psi_j||_{L^p_\mu(B_j)} &\leq&
||f||_{L^{pt'}_\mu(B_j)} ||\sqrt{Q}\nabla\psi_j||_{L_\mu^s(B_j)}
 \nonumber\\
&\leq& \overline{C} C_1(B_j)||(f,\nabla f)||_{\wa{p}{\Omega}}.
\eea
Combining this with (\ref{3.6-2}) gives
\bea ||\psi_j f||_{L^{p\sigma}_w(B_j)} &\leq& C(B_j)\big(1+\overline{C}
C_1(B_j)\big)||(f,\nabla f)||_{\wa{p}{\Omega}}.\nonumber
\eea
By (\ref{3.6-1}), for any $f\in {\cal H}$,
\bea \nonumber||f||_{L^{p\sigma}_w(\Omega')}&\leq& ||
(f,\nabla f)||_{\wa{p}{\Omega}}\disp\sum_j C(B_j)\big(1+ \overline{C}
C_1(B_j)\big) \nonumber\\ \nonumber &=& C(\Omega') ||(f,\nabla
f)||_{\wa{p}{\Omega}},
\eea
which completes the proof of Lemma \ref{localglobal}. \qed
\\

Theorem \ref{newappth2} follows from Lemma \ref{localglobal}
and Theorem \ref{generallocal}. We will sketch the proof, omitting
some familiar details. By choosing a sequence of compact sets
increasing to $\Omega$ and using a diagonalization argument, it is
enough to prove the conclusion for a fixed measurable $\Omega'$ with
compact closure $\overline{\Omega'}$ in $\Omega$. Fix such an
$\Omega'$ and select a bounded open $\Omega''$ with 
$\overline{\Omega'} \subset \Omega'' \subset \overline{\Omega''}
\subset \Omega.$  For ${\cal H}$ as in Theorem \ref{newappth2}, apply
Lemma \ref{localglobal} to the set $\Omega''$ to obtain
\bea\label{Omegaprimeprime}
||f||_{L^{p\sigma}_w(\Omega'')} \le C(\Omega'')\, ||(f,\nabla
  f)||_{W^{1,p}_{\nu,\mu}(\Omega,Q)},\quad f \in {\cal H}.
\eea
By assumption, $w<<\nu$, so (\ref{Omegaprimeprime}) extends to
$\overline{\cal H}$ in the form
\bea\label{Omegaprimeprimeext}
||f||_{L^{p\sigma}_w(\Omega'')} \le C(\Omega'')\,
  ||(f,\vec{g})||_{W^{1,p}_{\nu,\mu}(\Omega,Q)},\quad (f,\vec{g})\in
  \overline{\cal H}.
\eea

Let $\epsilon >0$. By hypothesis, ${\cal H}$ satisfies the Poincar\'e
estimate (\ref{poincare}) for balls $B_r(y)$ with $y \in
\overline{\Omega'}$ and $r<\delta(\epsilon, \Omega')$. Since the
Euclidean distance between $\overline{\Omega'}$ and $\partial
\Omega''$ is positive and we have assumed (\ref{eucled2}), we may also
assume by Remark \ref{various}(ii) that all such balls lie in the
larger set $\Omega''$. Next we claim that (\ref{poincare}) extends to
$\overline{\cal H}$, i.e.,
\bea\label{poincareext}
\left(\int_{B_r(y)} |f-f_{B_r(y),w}|^p dw\right)^{\frac{1}{p}} \le
\epsilon ||(f,\vec{g})||_{W^{1,p}_{\nu,\mu}(B_{c_0r}(y),Q)} \quad
\text{if } (f,\vec{g}) \in \overline{\cal H},
\eea
for the same class of balls $B_r(y)$. In fact, if $(f,\vec{g}) \in
\overline{\cal H}$ and $\{f^j\} \subset {\cal H}$ satisfies
$(f^j, \nabla f^j) \ra (f,\vec{g})$ in $W^{1,p}_{\nu,\mu}(\Omega,Q)$
norm, then there is a subsequence, still denoted $\{f^j\}$, with $f^j
\ra f$ a.e.-$\nu$ in $\Omega$, and so with $f^j \ra f$ a.e.-$w$ in
$\Omega$ since $w<<\nu$. By (\ref{Omegaprimeprime}), $\{f^j\}$ is
bounded in $L^{p\sigma}_w(\Omega'')$. Hence, since the balls in
(\ref{poincareext}) satisfy $B_r(y) \subset \Omega''$, we obtain
$f^j_{B_r(y),w} \ra f_{B_r(y),w}$ by our usual weak convergence
argument, and (\ref{poincareext}) follows by Fatou's lemma from its
analogue (\ref{poincare}) for the $(f^j, \nabla f^j)$.

Now let $\{(f_k,\vec{g_k})\} \subset \overline{\cal H}$ be bounded in
$W^{1,p}_{\nu,\mu}(\Omega,Q)$ norm and apply Theorem
\ref{generallocal} with $\mathcal{X}(\Omega) = L^p_\nu(\Omega) \times
{\cal L}^p_\mu(\Omega,Q)$ to the set ${\cal S}$ defined by
\bea
\nonumber {\cal S} = \left\{\big(f_k, (f_k,\vec{g_k})\big)\right\}_k,
\eea
and with $\{(E_\ell^\epsilon, F_\ell^\epsilon)\}_\ell$ chosen to be a
finite number of pairs $\{(B_r(y_\ell), B_{c_0r}(y_\ell)\}_\ell$ as in
(\ref{poincareext}), but now with $r$ fixed depending on $\epsilon$,
and with $\Omega' \subset \cup_\ell B_r(y_\ell)$. Such a finite choice
exists by (\ref{new2}) and the Heine-Borel theorem since
$\overline{\Omega'}$ is compact; cf. the proof of Lemma
\ref{help}. Since $\Omega'$ is completely covered by $\cup_\ell
E^\epsilon_\ell$, assumption (i) of Theorem \ref{generallocal} is
fulfilled. Moreover, the collection $\{F_\ell^\epsilon\}$ has bounded
overlaps uniformly in $\epsilon$ by the geometric doubling argument
used to prove Lemma \ref{help}.

Finally, (\ref{bddlocal}) follows from (\ref{Omegaprimeprimeext})
applied to the bounded sequence $\{(f_k,\vec{g_k})\}$ since
$\cup_{\ell, \epsilon} E_\ell^\epsilon \subset \Omega''$. Thus Theorem
\ref{generallocal} implies that there is a subsequence $\{f_{k_i}\}$
of $\{f_k\}$ and a function $f\in L^{p\sigma}_w(\Omega')$ such that
$f_{k_i} \ra f$ a.e.-$w$ in $\Omega'$ and in $L^q_w(\Omega')$ norm,
$1\le q < p\sigma$. This completes the proof of Theorem
\ref{newappth2}. \qed

For functions which are compactly supported in a fixed bounded
measurable $\Omega'$ with $\overline{\Omega'}\subset\Omega$, the proof
of Theorem \ref{newappth2} can be modified to yield compact embedding into
$L^q_w(\Omega')$ for the same $\Omega'$ without assuming
(\ref{eucled2}). Of course we always require (\ref{new2}). Given such 
$\Omega'$ and a set ${\cal H} \subset Lip_{Q,p,0}(\Omega')$, we may
view ${\cal H}$ as a subset of $Lip_{Q,p,0}(\Omega)$ simply by
extending functions in ${\cal H}$ to all of $\Omega$ as $0$ in
$\Omega\setminus\Omega'$. In this way, the proof of Theorem
\ref{newappth2} works without (\ref{eucled2}). For example, choosing
${\cal H} = Lip_{Q,p,0}(\Omega')$, we obtain 

\begin{thm}\label{last} Let the assumptions of \S 3.1 hold and
$w<<\nu$. Let $\Omega'$ be a bounded measurable set with
$\overline{\Omega'} \subset \Omega$. Fix $p\in [1,\infty)$ and suppose
the Poincar\'e property of order $p$ in Definition
\ref{poincaredef} holds for $Lip_{Q,p,0}(\Omega')$, with
$Lip_{Q,p,0}(\Omega')$ viewed as a subset of $Lip_{Q,p,0}(\Omega)$ using
extension by $0$, and suppose the local Sobolev property of order
$p$ in Definition \ref{localsobdef} holds. Assume the cutoff property
of some order $s\geq p\sigma'$ is true for $\mu$, with $\sigma$ as in
(\ref{sob}), and that (\ref{globpc}) holds for $Lip_{Q,p,0}(\Omega')$
with $t=s/p$. Then for every sequence $\{(f_k,\vec{g_k})\} \subset
W^{1,p}_{\nu,\mu,0}(\Omega',Q)$ which is bounded in
$W^{1,p}_{\nu,\mu}(\Omega',Q)$ norm, there is a subsequence
$\{f_{k_i}\}$ of $\{f_k\}$ and a function $f \in
L^{p\sigma}_w(\Omega')$ such that $f_{k_i} \ra f$ pointwise
a.e.-$w$ in $\Omega'$ and in $L^q_w(\Omega')$ norm, $1\le q
<p\sigma$.
\end{thm}

The full force of the local Sobolev estimate in Definition
\ref{localsobdef} is not needed to prove Theorem \ref{last}. In fact,
it is enough to assume that (\ref{sob}) holds only for balls centered
in the fixed compact set $\overline{\Omega'}$.

The proof of Theorem \ref{last} is like that of Theorem
\ref{newappth2}, working with the set $\Omega'$ that occurs in the
hypotheses of Theorem \ref{last}. However, now (\ref{localglobalsob})
in the conclusion of Lemma \ref{localglobal} (with ${\cal H} =
Lip_{Q,p,0}(\Omega')$) remains valid if $\Omega'$ is replaced on the
left side by $\Omega$ since every $f \in  Lip_{Q,p,0}(\Omega')$
vanishes on $\Omega\setminus\Omega'$. The resulting estimate serves as
a replacement for (\ref{Omegaprimeprime}), so it is not necessary to
demand that the $E_\ell^\epsilon$ are subsets of a compact set
$\overline{\Omega''} \subset \Omega$. Hence (\ref{eucled2}) is no
longer required. Finally, the Poincar\'e estimate extends as usual to
$W^{1,p}_{\nu,\mu,0}(\Omega', Q)$ (the closure of
$Lip_{Q,p,0}(\Omega'))$, and due to support considerations, 
the $E_\ell^\epsilon$ can be restricted to subsets of $\Omega'$ by
replacing $E_\ell^\epsilon$ by $E_\ell^\epsilon \cap \Omega'$; this
guarantees $w(E_\ell^\epsilon)< \infty$ since $w$ is locally finite by
hypothesis. 

Recalling the comments made immediately after Definition
\ref{cutoff} and in Remark \ref{extracondition}, we obtain a
useful special case of Theorem \ref{last}:

\begin{corollary}\label{lastcor}
Let the assumptions of \S 3.1 hold, $\Omega$ and $Q$ be bounded,
$w=\nu=\mu$ and (\ref{eucled2}) be true. Let $\Omega'$ be a measurable
set with $\overline{\Omega'} \subset \Omega$. Fix $p\in [1,\infty)$
and suppose the Poincar\'e property of order $p$ in Definition
\ref{poincaredef} holds for $Lip_{Q,p,0}(\Omega')$ and the local
Sobolev property of order $p$ in Definition \ref{localsobdef} holds.
Then for every $\{(f_k,\vec{g_k})\} \subset
W^{1,p}_{\nu,\mu,0}(\Omega',Q)$ which is bounded in
$W^{1,p}_{\nu,\mu}(\Omega,Q)$ norm, there is a subsequence
$\{f_{k_i}\}$ of $\{f_k\}$ and a function $f \in
L^{p\sigma}_w(\Omega')$ such that $f_{k_i} \ra f$ pointwise a.e.-$w$
in $\Omega'$ and in $L^q_w(\Omega')$ norm, $1\le q <p\sigma$.
\end{corollary}

In case $p=2$ and all measures are Lebesgue measure,  Corollary
\ref{lastcor} is used in \cite{R1} to show existence of weak solutions
to Dirichlet problems for some linear subelliptic equations. It is
also used in \cite{R2} to derive the global Sobolev inequality
\bea ||f||_{L^{2\sigma}(\Omega')} \leq C
\Big(\int_{\Omega'}|\sqrt{Q}\nabla f|^2dx\Big)^{1/2}
\eea
for open $\Omega'$ with $\overline{\Omega'}\subset\Omega$ from the
local estimate (\ref{sob2}).

\section{Precompact subsets of $L^N$ in a quasimetric space}
\setcounter{equation}{0}

\setcounter{theorem}{0}

\setcounter{equation}{0}

\setcounter{theorem}{0}
In this section, we will consider the situation of an open set
$\Omega$ in a topological space $X$ when $X$ is also endowed with a
quasimetric $d$. As there is no easy way to define Sobolev spaces on
general quasimetric spaces, this section concentrates on establishing
a simple criterion not directly related to Sobolev spaces
ensuring that bounded subsets of $L^N_w(\Omega)$ are precompact in
$L^q_w(\Omega)$ when $1\le q<N\le \infty$.

We begin by further describing the setting for our result. The
topology on $X$ is expressed in terms of a fixed collection
${\cal T}$ of subsets of $X$ which may not be related to the
quasimetric $d$. Thus when we say that a set ${\cal O}\subset X$ is
\emph{open}, we mean that ${\cal O}\in {\cal T}$. Given an open
$\Omega$, we will assume each of the following:   
\bea
\nonumber &&(i)\;\;\; \mbox{$\forall x\in X$ and $r>0$, the $d$-ball
  $B_r(x) = \{y\in X\;:\; d(x,y)<r\}$ is a  Borel set;}\\
\nonumber  &&(ii)\;\;\mbox{$\forall x\in X$ and $r>0$, there is an open
  set ${\cal O}$ so that $x\in{\cal  O}\subset B_r(x)$;}\\
\nonumber &&(iii)\;\;\mbox{if $X\not = \Omega$, then $\forall x\in
  \Omega$, $d(x,\Omega^c)=\inf\{ d(x,y): y\in \Omega^c\} >0$.}
\eea

Property $(ii)$ serves as a substitute for (\ref{new2}).

Unlike the situation in \S 3, $d$-balls centered in $\Omega$ may not
be subsets of $\Omega$ unless $X= \Omega$. However, we note the
following fact. 
\begin{remark}\label{remark1}
Properties $(ii)$ and $(iii)$ guarantee that for any compact set
$K\subset \Omega$, there exists $\vep(K)>0$ such that $B_r(x)\subset
\Omega$ if $x\in K$ and $r < \vep(K)$. In fact, first note that for
any $x\in \Omega$, $(iii)$ implies that the $d$-ball $B(x)$ with center
$x$ and radius $r_x = d(x,\Omega^c)/(2\kappa)$ lies in $\Omega$. If
$K$ is a compact set in $\Omega$, (ii) shows that $K$ can be covered
by a finite number of such balls $\{B(x_i)\}$. With $\vep(K)$ chosen
to be a suitably small multiple (depending on $\kappa$) of $min\,
\{r_{x_i}\}$, the remark then follows easily from the swallowing property
of $d$-balls. 
\end{remark}

Further, we assume that $(\Omega,d)$ satisfies the local geometric
doubling condition in Definition \ref{doublingdef}, i.e., for each
compact set $K\subset \Omega$, there exists $\delta'(K)>0$ such
that for all $x\in K$ and all $0<r'<r < \delta'(K)$, the number of
disjoint $d$-balls of common radius $r'$ contained in $B_{r}(x)$ is at
most a constant ${\cal C}_{r/r'}$ depending on $r/r'$ but not on $K$. 
We will choose $\delta'(K)\le \vep(K)$ in the above.  

With this framework in force, we now state the main result of the section.

\begin{thm}\label{mainmetric}
Let $\Omega\subset X$ be as above, and let $w$ be a finite Borel
measure on $\Omega$ such that given any $\epsilon >0$, there is a
compact set $K \subset \Omega$ with $w(\Omega\setminus K)<
\epsilon$. Let $1\leq p<\infty$ and $1<N\leq \infty$, and suppose
$\CS\subset L^N_w(\Omega)$ has the property that for any compact set
$K\subset \Omega$, there exists $\delta_K>0$ such that
\begin{equation}\label{pointype}
\|f-f_{B,w}\|_{L^p_w(B)}\le b(f,B) \ \mbox{ if $f\in \CS$ and
$B=B_r(x)$, $x\in K$, $0<r<\delta_K$},
\end{equation}
where $b(f,B)$ is a nonnegative ball set function.  Further, suppose
there is a constant $c_0\ge 1$ so that for every $\ep>0$ and
every compact set $K\subset \Omega$, there exists
$\tilde{\delta}_{\ep,K}>0$ such that
\begin{equation}\label{overlap}
\sum_{B\in \F} b(f,B)^p \le \ep^p \ \mbox{ for all $f\in \CS$ }
\end{equation}
for every finite family $\F=\{B\}$ of $d$-balls centered in $K$ with
common radius less than $\tilde{\delta}_{\ep,K}$ for which $\{c_0B\}$
is a pairwise disjoint family of subsets of $\Omega$. Then any
sequence in $\CS$ that is bounded in $L^N_w(\Omega)$ has a subsequence
that converges in $L^q_w(\Omega)$ for $1\le q<N$ to a function in
$L^N_w(\Omega)$. 
\end{thm}

{\bf Proof. } Let $\ep>0$ and choose a compact set $K\subset \Omega$ with
$w(\Omega\setminus K)<\ep$. Next, for $c_0\ge 1$, as in the proof of
Lemma \ref{help} there is a positive constant $r=r(\ep,K,c_0)<
\min\{\delta_{K},\tilde{\delta}_{\ep,K},\delta'(K),\vep(K)/(\gamma
c_0)\}$ (see (\ref{pointype}),(\ref{overlap}), Definition
\ref{doublingdef} and Remark \ref{remark1}), where $\gamma=\k+2\k^2$
with $\k$ as in (\ref{tri}), and a finite family $\{B_r(y_k)\}_k$ 
of $d$-balls centered in $K$ satisfying $K\subset \cup_k
B_r(y_k)$ and whose dilates $\{B_{c_0r}(y_k)\}_k$ lie
in $\Omega$ and have the bounded intercept property (with intercept
constant $M$ independent of $\ep$).  Since $\{B_{c_0r}(y_k)\}_k$ has bounded intercepts with
bound $M$, it can be written as the union of at most $M$ families of
disjoint $d$-balls; see e.g. the proof of \cite[Lemma 2.5]{CW1}. By
(\ref{overlap}), we conclude that
$$
\sum_{k} b(f,B_r(y_k))^p\le M\ep^p.
$$
Theorem \ref{mainmetric} then follows immediately from Theorem
\ref{absgeneral}; see also Remark 1.3(1). \qed
 
As an application of Theorem \ref{mainmetric} we present a version of
\cite[Theorem 8.1]{HK2} in the case $p\geq 1$. Our version improves
the one in \cite{HK2} by allowing two different measures and by
relaxing the assumptions made about embedding and
doubling. Furthermore, while the analogue in \cite{HK2} of our
(\ref{pwitha}) uses only the $L^1_w(B)$ norm on the left side, it
automatically self-improves to the $L^p_w(B)$ norm due to the 
doubling assumption, with a further fixed enlargement of the ball
$c_0B$ on the right side; see e.g. \cite[Theorem 5.1]{HK2}.

\begin{corollary} \label{betterthanHK} Let $X,d,\Omega,w$ be as above,
  and let $\mu$ be a Borel measure on $\Omega$. Fix   $1\leq
  p<\infty$, $1<N\leq \infty$ and $c_0\geq 1$. Consider a 
  sequence of pairs $\{(f_i,g_i)\}\subset L^N_w(\Omega)\times
  L^p_\mu(\Omega)$ such that for any compact set $K\subset \Omega$,
  there exists $\bar{\delta}_K>0$ with 
\bea\label{pwitha} ||f_i-(f_i)_{B,w}||_{L^p_w(B)} \leq
  a_*(B)||g_i||_{L^p_\mu(c_0B)} 
\eea
for all $i$ and all $d$-balls $B$ centered in $K$ with $c_0B
  \subset\Omega$ and $r(B)<\bar{\delta}_K$, where $a_*(B)$ is a
  non-negative ball set function satisfying 
\bea\label{a_*}
\disp\lim_{r\ra 0}\Big\{ \sup_{y\in K}a_*(B_r(y))\Big\}=0.
\eea
Then if $\{f_i\}$ and $\{g_i\}$ are bounded in $L^N_w(\Omega)$ and
$L^p_\mu(\Omega)$ respectively, $\{f_i\}$ has a subsequence
converging in $L^q_w(\Omega)$ for $1\leq q<N$ to a function belonging
to $L^N_w(\Omega)$.
\end{corollary}

{\bf Proof.}  Given $\ep>0$ and compact set $K\subset \Omega$, use
(\ref{a_*}) to choose $r_0>0$ so that $a_*(B_r) < \epsilon/\beta$ for
any $d$-ball $B_r$ centered in $K$ with $r<r_0$, where $\beta =
\sup_i||g_i||_{L^p_\mu(\Omega)}<\infty$.  In Theorem \ref{mainmetric},
choose $\CS = \{f_i\}$, $\delta_K = \overline{\delta}_K$, 
$b(f_i,B) = a_*(B)||g_i||_{L^p_\mu(c_0B)}$ and
$$
\tilde{\delta}_{\ep,K} = \min\{\overline{\delta}_K, \delta'(K)
,r_0,\vep(K)/c_0\}.
$$
If $B$ is a $d$-ball with center in $K$ and $r(B)
<\tilde{\delta}_{\ep,K}$, then $c_0B\subset \Omega$. Hence, 
\bea\label{4.10} \sum_{B\in{\cal F}}\big(a_*(B)||g_i||_{L^p_\mu(c_0B)}
\big)^p\leq \ep^p||g_i||_{L^p_\mu(\Omega)}^p/\beta^p \le
\ep^p\nonumber 
\eea
for every ${\cal F}$ as in Theorem \ref{mainmetric}. The conclusion
now follows from Theorem \ref{mainmetric}.\qed

\begin{remark}
\begin{enumerate}
\item  The $g_i$ in (\ref{pwitha}) are usually the modulus of a
fixed derivative of the corresponding $f_i$, such as $|\nabla f_i|$
when $X$ is a Riemannian manifold. More generally, $g_i$ may be the
upper gradient of $f_i$ (see \cite{Hei} for  the 
definition).
\item Theorem \ref{mainmetric} can also be used to obtain an extension
of Theorem 2.3 to $s$-John domains in quasimetric spaces; see
\cite[Theorem 1.6]{CW2}.
\end{enumerate}
\end{remark}

\section{Appendix}
\setcounter{equation}{0}
\setcounter{theorem}{0}

Here we briefly consider analogues of Theorem \ref{appth1}, Corollary
\ref{appcor1} and Theorem \ref{newappth2} without assuming $w<<\nu$,
but adding the assumption that ${\cal H}$ is linear. In this case,
(\ref{strongsob}) can be extended by continuity to obtain a bounded
linear map from $\overline{\cal H}$ into
$L^{p\sigma}_w(\Omega)$. Here, as always, $\overline{\cal H}$ denotes 
the closure of $\{(f,\nabla f): f\in {\cal H}\}$ in
$W^{1,p}_{\nu,\mu}(\Omega, Q)$. However, when $w<<\nu$ fails, there is
no natural way to obtain the extension for every $(f,\vec{g}) \in
\overline{\cal H}$ keeping the same $f$ on the left side. In fact, let
$(f,\vec{g}) \in \overline{{\cal H}}$ and choose $\{f_j\}\subset {\cal
H}$ with $(f_j, \nabla f_j)\to (f,\vec{g})$ in
$W^{1,p}_{\nu,\mu}(\Omega,Q)$. Linearity of ${\cal H}$ allows us to
apply (\ref{strongsob}) to differences of the $f_j$ and conclude that
$\{f_j\}$ is a Cauchy sequence in $L^{p\sigma}_w(\Omega)$. Therefore
$f_j \to f^*$ in $L^{p\sigma}_w(\Omega)$ for some $f^* \in
L^{p\sigma}_w(\Omega)$, and 
\[
||f^*||_{L^{p\sigma}_w(\Omega)} \le C ||(f,\vec{g})||_{\wa{p}{\Omega}}
\quad\text{if $(f,\vec{g})\in \overline{\cal H}$}.
\]
The function $f^*$ is determined by $(f, \vec{g})$, i.e., $f^*$
is independent of the particular sequence $\{f_j\}\subset
{\cal H}$ above. Indeed, if $\{\tilde{f_j}\}$ is another sequence in
${\cal H}$ with $(\tilde{f_j},\nabla\tilde{f_j})
\ra (f, \vec{g})$ in $W^{1,p}_{\nu,\mu}(\Omega, Q)$, and if
$\tilde{f_j} \ra \tilde{f^*}$ in $L^{p\sigma}_w(\Omega)$, then by
(\ref{strongsob}) and linearity of ${\cal H}$,
\[
||\tilde{f_j}- f_j||_{L^{p\sigma}_w(\Omega)} \le C
||(\tilde{f_j}-f_j,\nabla \tilde{f_j} - \nabla f_j)||_{W^{1,p}_{\nu,
\mu}(\Omega, Q)} \ra 0.
\]
Consequently $||\tilde{f^*}- f^*||_{L^{p\sigma}_w(\Omega)} =0$. Thus
$(f, \vec{g})$ determines $f^*$ uniquely as an element of
$L^{p\sigma}_w(\Omega)$. Define a mapping
\bea\label{map}
T: \overline{\cal H} \ra L^{p\sigma}_w(\Omega)\quad\text{by setting
$T(f, \vec{g}) = f^*$}.
\eea
Note that $\overline{\cal H}$ is a linear set in
$W^{1,p}_{\nu,\mu}(\Omega, Q)$ since ${\cal H}$ is linear, and that
$T$ is a bounded linear map from $\overline{\cal H}$ into
$L^{p\sigma}_w(\Omega)$.  Also note that $T$
satisfies $T(f,\nabla f) = f$ when restricted to those $(f, \nabla f)$
with $f\in {\cal H}$. Furthermore, if $w <<\nu$ then $T(f,\vec{g}) =f$ for
all $(f,\vec{g}) \in \overline{\cal H}$, i.e., $f^* = f$ a.e.-$w$ for all
$(f,\vec{g}) \in \overline{\cal H}$.  This follows since $f_j \ra f$
in $L^p_{\nu}(\Omega)$ norm and $f_j \ra f^*$ in
$L^{p\sigma}_w(\Omega)$ norm. In this appendix, where it is not
assumed that $w<<\nu$, $f^*$ plays a main role. One can find a
function $h$ such that $h=f^*$ a.e.-$w$ and $h=f$ a.e.-$\nu$, but as
this fact is not needed, we omit its proof.

An analogue of Theorem \ref{appth1} is given in the next result. 

\begin{thm} \label{appth1*}   Let all the assumptions of Theorem
\ref{appth1} hold except that now the set ${\cal H}$ is linear and we
do not assume $w<<\nu$. Then the map $T: \overline{{\cal H}} \ra
L^q_w(\Omega)$ defined in (\ref{map}) is compact if $1\leq q <
p\sigma$. Equivalently, if $\{(f_k,\vec{g_k})\}$ is a sequence in
$\overline{{\cal H}}$ with $\sup_k
||(f_k,\vec{g_k})||_{W^{1,p}_{\nu,\mu}(\Omega,Q)} < \infty$, then
$\{f^*_k\}$ has a subsequence which converges in $L^q_w(\Omega)$ for
$1 \le q< p\sigma$, where $f^*_k= T(f_k,\vec{g_k})$. Moreover, the
limit of the subsequence belongs to $L^{p\sigma}_w(\Omega)$.
\end{thm}
{\bf Proof:}  Let ${\cal H}$ satisfy the hypothesis of the theorem and
let $\{(f_k,\vec{g_k})\}\subset \overline{\cal H}$ be bounded in
$W^{1,p}_{\nu,\mu}(\Omega, Q)$.  For each $k$, choose $h_k\in {\cal
H}$ so that 
\bea\label{nearby} ||(f_k,\vec{g_k})-(h_k,\nabla
h_k)||_{\wa{p}{\Omega}} \leq 2^{-k}.
\eea
Set ${\cal H}_1= \{h_k\}_k \subset {\cal H}$.  Then
$\{(h_k, \nabla h_k): h_k \in {\cal H}_1\}$ is bounded in
$\wa{p}{\Omega}$. Further, (\ref{strongsob}) implies a version of
(\ref{sup}), namely
\[
\sup_{f \in {\cal H}_1}
\left\{||f||_{L^{p\sigma}_w(\Omega)} + ||(f,\nabla
f)||_{W^{1,p}_{\nu,\mu}(\Omega,Q)} \right\}  < \infty.
\]
Theorem \ref{simpleversion} now applies to ${\cal H}_1$ with
$N=p\sigma$ and gives that any sequence in $\hat{{\cal H}_1}$ has a
subsequence which converges in $L^q_w(\Omega)$ norm for $1\leq q<
p\sigma$ to a function belonging to $L^{p\sigma}_w(\Omega)$.
The sequence $\{h_k\}$ lies in $\hat{{\cal H}_1}$, as is
easily seen by considering, for each fixed $k$, the constant sequence
$\{f^j\}$ defined by $f^j = h_k$ for all $j$. We conclude that
$\{h_k\}$ has a subsequence $\{h_{k_l}\}$ converging in
$L^q_w(\Omega)$ norm for $1\le q <p\sigma$ to a function $h
\in L^{p\sigma}_w(\Omega)$. By linearity and boundedness of $T$ from
$\overline{\cal H}$ to $L^{p\sigma}_w(\Omega)$ together with
(\ref{nearby}), we have (writing $f_k^* = T(f_k, \vec{g_k})$)
\[
||f_k^* -h_k||_{L^{p\sigma}_w(\Omega)} = ||T(f_k,\vec{g_k}) - T(h_k,
  \nabla h_k)||_{L^{p\sigma}_w(\Omega)} \le C 2^{-k} \rightarrow 0.
\]
Restricting $k$ to $\{k_l\}$ and using $w(\Omega)<\infty$, we conclude
that $\{f^*_{k_l}\}$ also converges to $h$ in $L^q_w(\Omega)$ for
$1\le q<p\sigma$, which completes the proof. \qed  

Setting ${\cal H}= Lip_{Q,p}(\Omega)$ in Theorem \ref{appth1*} gives
an analogue of Corollary \ref{appcor1}: 

\begin{corollary}\label{appcor1*} Let the hypotheses of
Theorem \ref{appth1*} hold for ${\cal H} = Lip_{Q,p}(\Omega)$.
Then the map $T$ defined by (\ref{map}) is a compact map of
$\wa{p}{\Omega}$ into $L^q_w(\Omega)$ for $1\leq q < p\sigma$,
i.e., if $\{(f_k,\vec{g_k})\} \subset \wa{p}{\Omega}$ and
$\sup_k ||(f_k,\vec{g_k})||_{W^{1,p}_{\nu,\mu}(\Omega,Q)} < \infty$,
then $\{ f^*_k\}$ has a subsequence which converges in
$L^q_w(\Omega)$ for $1 \le q< p\sigma$, where $f^*_k =
T(f_k,\vec{g_k})$. Moreover, the limit of the subsequence belongs to
$L^{p\sigma}_w(\Omega)$.
\end{corollary}

Theorem \ref{newappth2} also has an analogue without assuming $w<<\nu$
provided ${\cal H}$ is linear, and in this instance (\ref{strongsob})
is not required: the subsequence $\{f_{k_i}\}$ of $\{f_k\}$ in
the conclusion is then replaced by a subsequence of $\{f_k^*\}$, where
$f_k^*$ is constructed as above but now using bounded measurable
$\Omega'$ whose closures increase to $\Omega$. Now $f^*$ arises when
(\ref{Omegaprimeprime}) is extended to $\overline{\cal H}$, namely,
instead of (\ref{Omegaprimeprimeext}), we obtain
\bea
\nonumber
||f^*||_{L^{p\sigma}_w(\Omega'')} \le C(\Omega'')\,
  ||(f,\vec{g})||_{W^{1,p}_{\nu,\mu}(\Omega,Q)}\quad \text{if
}(f,\vec{g})\in   \overline{\cal H}
\eea
where $f^*$ is constructed for a pair $(f,\vec{g})\in \overline{\cal
H}$ by using linearity of ${\cal H}$ and (\ref{Omegaprimeprime}) for a
particular $(\Omega', \Omega'')$. It is easy to see that $f^* \in
L^{p\sigma}_{w, loc}(\Omega)$ by letting $\Omega' \nearrow
\Omega$. The Poincar\'e inequality analogous to (\ref{poincareext}) is
\bea
\nonumber \left(\int_{B_r(y)} |f^*-f^*_{B_r(y),w}|^p dw
\right)^\frac{1}{p} \le
\epsilon ||(f,\vec{g})||_{W^{1,p}_{\nu,\mu}(B_{c_0r}(y),Q)} \quad
\text{if } (f,\vec{g}) \in \overline{\cal H},
\eea
obtained by extending (\ref{poincare}) from ${\cal H}$ to
$\overline{\cal H}$. Further details are omitted.\\

\singlespace
{\footnotesize{

\bigskip

\noindent Department of Mathematics\\
National University of Singapore\\
10, Lower Kent Ridge Road\\
Singapore 119076\\
e-mail: matcsk@nus.edu.sg\\

\noindent Department of Mathematics, Physics and Geology\\
Cape Breton University\\
Sydney, NS B1P6L2\\
e-mail: scott.rodney@gmail.com\\

\noindent Department of Mathematics\\
Rutgers University\\
Piscataway, NJ 08854\\
e-mail: wheeden@math.rutgers.edu

}

\end{document}